\documentclass[pdflatex,sn-mathphys-ay]{sn-jnl}

\usepackage{graphicx}%
\usepackage{multirow}%
\usepackage{amsmath,amssymb,amsfonts}%
\usepackage{amsthm}%
\usepackage[title]{appendix}%
\usepackage{textcomp}%
\usepackage{manyfoot}%
\usepackage{booktabs}%

\usepackage{pgfplots}
\pgfplotsset{compat=newest}
\pgfplotsset{plot coordinates/math parser=false}
\newlength\figureheight
\newlength\figurewidth

\usepackage[labelformat=simple,caption=false]{subfig}
\usepackage{xcolor}
\usepackage{hyperref}        
\usepackage{tikz,tikz-3dplot}
\usetikzlibrary{calc,3d,shapes,arrows.meta,perspective,positioning,fit}
\hypersetup{colorlinks, breaklinks,
			linkcolor=blue, 
			citecolor=black, 
			urlcolor=blue,
            pdfstartview=FitV}
\usepackage[capitalise, noabbrev]{cleveref}

\newcommand{\vect}[1]{\ensuremath{\boldsymbol{#1}}}

\newcommand{\rt}[2]{\ensuremath{#1{:}#2}}

\begin{document}

\title[Cislunar State and Uncertainty Propagation via the Modified Generalized Equinoctial Orbital Elements]{Cislunar State and Uncertainty Propagation via the Modified Generalized Equinoctial Orbital Elements}

\author*[1]{\fnm{Maaninee} \sur{Gupta}}\email{maanineegupta@gmail.com}

\author[2]{\fnm{Kyle} J. \sur{DeMars}}\email{kdemars@purdue.edu}

\affil[1,2]{\orgdiv{School of Aeronautics and Astronautics}, \orgname{Purdue University}, \orgaddress{\city{West Lafayette}, \state{IN} \postcode{47906}, \country{USA}}}

\abstract{The complex cislunar dynamical environment poses challenges for spacecraft navigation and Space Domain Awareness (SDA) operations, where the knowledge of current and future spacecraft states is essential. Conventional Gaussian-based approaches for SDA degrade under the nonlinearities that manifest in this regime. To accurately model the underlying dynamics and characterize uncertainty, this work explores the Modified Generalized Equinoctial Orbital Elements under high-fidelity propagation for cislunar applications. The Henze--Zirkler test for multivariate normality is leveraged to evaluate uncertainty evolution across a range of orbits, demonstrating improved preservation of Gaussian behavior in cislunar space.}


\maketitle

\section{Introduction}\label{sec1}

With the recent interest in missions to the vicinity of the Moon, cislunar space is poised to become the domain that sustains humanity's presence beyond the Earth. With this renewed interest, Space Domain Awareness (SDA), that has conventionally applied to the sub-geosynchronous orbit regime, will be necessary to support cooperative and safe operations in cislunar space as well. In contrast to SDA near the vicinity of the Earth, SDA across the broad volume of cislunar space presents unique challenges. Operationally, the vastness of this region, which is defined as the spherical volume contained within the orbit of the Moon around the Earth, introduces significant difficulties in observing, tracking, and maintaining custody of objects from Earth-based assets.

In addition to these operational challenges, the cislunar environment is characterized by highly nonlinear and chaotic dynamics arising from the non-negligible impact of lunar gravity, coupled with other perturbing effects. These complex dynamics complicate trajectory design and spacecraft navigation, where the accurate knowledge of both current and future spacecraft states is sought. Thus, realistic modeling of cislunar motion requires dynamical formulations that capture the influence of both the Earth and the Moon, along with other relevant perturbations, for long-term propagation and prediction. 

A critical component of cislunar SDA is the characterization and propagation of uncertainty associated with object state estimates. Accurate uncertainty representation is necessary for catalog maintenance, conjunction assessment, and maneuver detection \citep{Poore2016}. One factor that introduces inaccuracies in representing uncertainty over time is the assumption of Gaussian behavior, a simplification that degrades rapidly under the strong nonlinearities of this regime. As uncertainty distributions evolve, departures from Gaussianity introduce errors into estimation, prediction, and risk assessment frameworks that rely on Gaussian assumptions. Alternative state representations, particularly through the use of orbital element coordinates rather than conventional Cartesian position–velocity states, can ameliorate these effects. In cislunar space, however, careful treatment of the underlying dynamics remains essential to prevent divergence between predicted and true downstream uncertainty.

While osculating Keplerian orbital elements may be determined at discrete time steps for cislunar trajectories, their direct use for propagation requires explicit incorporation of the relevant perturbing accelerations into the equations of motion, an approach that is not commonly employed. As such, element sets that embedded orbital perturbations directly into their formulation have been explored. In particular, Broucke and Cefola showed that Equinoctial Orbital Elements (EqOEs) admit variational equations in Lagrange form that facilitate the inclusion of general perturbing accelerations directly within the element dynamics \citep{Broucke1972,Cefola1972}. An offshoot of the classical Equinoctial Orbital Elements includes the Modified Equinoctial Orbital Elements, wherein the semi-major axis and mean longitude from the classical set are replaced by the semi-latus rectum and the true longitude \cite{Walker1983}. \cite{Horwood2011} employ the Alternate Equinoctial Orbital Elements (AEqOEs) for their preservation of Gaussianity of the initial state uncertainty through propagation under Keplerian motion. Other developments of the EqOE framework, such as the oblate spheroidal EqOEs \citep{Biria2019} and $J_2$EqOEs \citep{Aristoff2021}, provide options for element sets that accommodate $J_2$ effects near the vicinity of the Earth. For modeling cislunar astrodynamics, however, there is a need for an orbital element set that encapsulates the lunar gravity. 

The Generalized Equinoctial Orbital Elements (GEqOEs) provide an avenue for a low-complexity methodology to incorporate both conservative and non-conservative perturbations that manifest in cislunar space. Introduced by Ba\`u et al., the GEqOE set has been successfully leveraged for state propagation of near-Earth orbits subject to third-body perturbations and $J_2$ effects \citep{Bau2021}. The GEqOE set is characterized by the total energy of the system and, as such, may be utilized to directly embed any perturbations arising from conservative and non-conservative forces directly into the orbital elements. For near-Earth applications, \cite{HernandoAyuso2023} and \cite{McGee2023} leverage GEqOEs for accurate state and uncertainty propagation under the addition of $J_2$ and third-body gravity effects. More recently, Gupta and DeMars have applied the GEqOEs for capturing three-body motion in cislunar space, demonstrating improved preservation of Gaussian behavior for uncertainty propagated along various cislunar orbits \citep{GuptaDeMars2025Astro,GuptaDeMars2025Unc}. While the GEqOE set demonstrates success in propagating various cislunar periodic orbits, its validity is restricted to trajectories characterized by a negative total energy. The Modified Generalized Equinoctial Orbital Elements (M-GEqOEs) supply the modifications necessary to overcome energy limitations and improve the robustness of modeling trajectories in the cislunar domain \citep{GuptaDeMars2025ESA}. 

Building on these advances, the current work extends the application of M-GEqOEs to high-fidelity cislunar dynamics with full Earth–Moon–Sun perturbations to evaluate their effectiveness for uncertainty characterization. Trajectories influenced by the combined gravitational effects of the Earth, Moon, and Sun are propagated directly in the generalized coordinates and validated against a Cartesian $N$-body ephemeris model. A range of orbits spanning the cislunar domain is examined to demonstrate effectiveness across the variety of dynamical structures that persist in this regime. In addition to assessing state propagation accuracy, this work focuses on the evolution of uncertainty under high-fidelity dynamics. Monte Carlo simulations are employed in both M-GEqOE and Cartesian coordinates, and the Henze-–Zirkler test for multivariate normality is used to quantify the preservation of Gaussianity. Through this analysis, the work provides insight into the role of coordinate choice in uncertainty propagation and its implications for cislunar SDA.

\section{The Modified Generalized Equinoctial Orbital Elements}
The Modified Generalized Equinoctial Orbital Elements (M-GEqOEs), derived from the Generalized Equinoctial Orbital Elements (GEqOEs), are introduced. The mathematical formulation of these elements, together with their time derivatives and the required coordinate transformations, is detailed in a general context.

\subsection{Generalized Orbital Motion}
Consider an object of mass $m$ orbiting the Earth or the Moon. In the associated body-centered inertial frame, the perturbed two-body motion of the object is represented as
\begin{align}\label{eq:pert2b}
    \ddot{\vect{r}} + \frac{\mu_{_{CB}} \vect{r}}{r^3} = \vect{a}_p(\vect{r}, \dot{\vect{r}}, t) \ ,
\end{align}
where $\vect{r}$ represents the inertial position vector of the object, $\dot{\vect{r}}$ is the inertial velocity, and $\ddot{\vect{r}}$ is the inertial acceleration. The vector $\vect{a}_p$ represents any perturbing accelerations acting on the object. The variables $r$ and $t$ denote the magnitude of the object position and time respectively. The parameter $\mu_{_{CB}}$ represents the gravitational parameter for the central body (either the Earth or the Moon), where the mass of the object is assumed to be infinitesimal relative to the mass of the central body. 

The origin of the inertial reference frame, denoted $O$, is centered on the celestial body and is defined as
\begin{align}\label{eq:inerframe}
    \Sigma = \{O; \ \vect{e}_x, \ \vect{e}_y, \ \vect{e}_z \} \ .
\end{align}
In \cref{eq:inerframe}, $\vect{e}_x = [1, \ 0, \ 0]^T$, $\vect{e}_y = [0, \ 1, \ 0]^T$, and $\vect{e}_z = [0, \ 0, \ 1]^T$. Additionally, the orbital reference frame, denoted $\Sigma_{or}$, is defined by the orthonormal basis,
\begin{align}\label{eq:orbrefframe}
    \Sigma_{or} = \{O; \ \vect{e}_r, \ \vect{e}_f, \ \vect{e}_h \} \ ,
\end{align}
where $\vect{e}_r$ points along the object position vector, $\vect{e}_h$ is directed along the angular momentum vector, and $\vect{e}_f$ completes the dextral orthonormal triad, 
\begin{align}
    \vect{e}_r = \frac{\vect{r}}{r}, \ \ \ \vect{e}_f = \vect{e}_h \times \vect{e}_r, \ \ \ \vect{e}_h = \frac{\vect{r} \times \dot{\vect{r}}}{|\vect{r} \times \dot{\vect{r}}|} = \frac{\vect{h}}{h} \ .
\end{align}
The perturbing accelerations in \cref{eq:pert2b} arise from a total perturbing force, denoted $\vect{F}$, that can be written as a sum of forces that are derived from a potential energy, such as the gravitational acceleration from additional celestial bodies, and forces that are not derived from a potential, such as atmospheric drag and solar radiation pressure. Mathematically, the total perturbing force is then represented as
\begin{align}\label{eq:pertforces}
    \vect{F} = \vect{P} - \nabla U \ ,
\end{align}
where $\vect{P}$ represents the perturbing forces that do not arise from a potential, and $- \nabla U$ models the contribution of forces derived from a potential energy, $U$. The total orbital energy of the object, $\mathcal{E}$, is a sum of its Keplerian energy, $\mathcal{E}_K$, and this potential energy, such that
\begin{align}\label{eq:totE}
    \mathcal{E} & = \mathcal{E}_K + U \ .
\end{align}
Expressing the object velocity vector, $\dot{\vect{r}}$, in the orbital frame as
\begin{align}
    \dot{\vect{r}} = \dot{r} \vect{e}_r + \frac{h}{r} \vect{e}_h \ ,
\end{align}
the total orbital energy may be represented as
\begin{align}
    \mathcal{E} = \frac{\dot{r}^2}{2} + \frac{h^2}{2r^2} - \frac{\mu_{_{C}}}{r} + U \ ,
\end{align}
where $\dot{r}$ denotes the radial component of the velocity, and $h$ is the magnitude of the angular momentum of the object. Then, the effective potential energy, denoted $U_{\mathrm{eff}}$, is defined as 
\begin{align}\label{eq:Ueff}
    U_{\mathrm{eff}} = \frac{h^2}{2r^2} + U \ .
\end{align}
Thus, the total energy of the object is
\begin{align}
    \mathcal{E} = \frac{\dot{r}^2}{2} - \frac{\mu_{_{C}}}{r} + U_{\mathrm{eff}} \ .
\end{align}
Using the effective potential energy, the \textit{generalized} angular momentum of the object, denoted $\Tilde{h}$, is introduced,
\begin{align}\label{eq:genangmom}
    \Tilde{h} = \sqrt{2 r^2 U_{\mathrm{eff}}}
\end{align}
which defines the \textit{generalized} velocity vector
\begin{align}
    \tilde{\vect{v}} = \dot{r} \vect{e}_r + \frac{\tilde{h}}{r} \vect{e}_h \ .
\end{align}
The generalized eccentricity vector, $\tilde{\vect{e}}$, orients the ellipse in the orbital plane and is given as
\begin{align}\label{eq:geneccvec}
    \mu_{_{C}} \tilde{\vect{e}} = \tilde{\vect{v}} \times (\vect{r} \times \tilde{\vect{v}} ) - \mu_{_{C}} \vect{e}_r \ .
\end{align}
Parameters denoted by tildes represent generalized quantities that are defined by the total orbital energy, $\mathcal{E}$, and are, thus, embedded with conservative perturbations acting on the object. Together, the vectors $(\vect{r}, \tilde{\vect{v}})$ characterize a non-osculating ellipse with one focus located at the center of the central gravitational body and perturbed by additional forces. Although the preceding formulation is provided in the context of orbits centered on the Earth or the Moon, the equations are applicable to orbits centered on other celestial bodies as well. 

\subsection{Modified Generalized Equinoctial Orbital Elements}
The set of Modified Generalized Equinoctial Orbital Elements (M-GEqOEs) is defined as
\begin{align}
    \{\tilde{p}, \ \ p_1, \ \ p_2, \ \ q_1, \ \ q_2, \ \ L \} \ .
\end{align}
The first element in the M-GEqOE set is the generalized semi-latus rectum, $\tilde{p}$, determined as
\begin{align}\label{eq:gensemilatrec}
    \tilde{p} = \frac{\tilde{h}^2}{\mu_{_{C}}} \ ,
\end{align}
where $\tilde{h}$ represents the generalized angular momentum, defined in \cref{eq:genangmom}, and $\mu_{_{C}}$ is the gravitational parameter associated with the central gravitational body. The sixth element in the M-GEqOE set is the classical true longitude, $L$, which represents the time-varying or ``fast'' variable along the orbit. The true longitude is determined as
\begin{align}{\label{eq:truelong}}
    L = \omega + \Omega + \theta \ ,
\end{align}
where $\theta$ is the classical true anomaly, and $\omega$ and $\Omega$ denote the argument of periapsis and the right ascension of ascending node, respectively. The second and third elements parameterize the eccentricity vector and are defined as
\begin{align}
    p_1 & = \tilde{e} \sin{\Psi} \label{eq:p1} \\ 
    p_2 & = \tilde{e} \cos{\Psi} \label{eq:p2} \ .
\end{align}
Here, the angle $\Psi$ denotes the generalized longitude of periapsis and is determined as
\begin{align}
    \Psi = L - \tilde{\theta} = \omega + \Omega + \theta - \tilde{\theta} \ ,
\end{align}
where $\tilde{\theta}$ is the generalized true anomaly. In the absence of orbital perturbations ($U = 0$), the generalized and classical true anomalies are equal, and $\Psi = \omega + \Omega$. Finally, the elements $q_1$ and $q_2$ orient the equinoctial reference frame relative to the inertial reference frame. These elements are functions of the classical inclination, $i$, and the classical right ascension of the ascending node, $\Omega$, expressed as
\begin{align}
    q_1 & = \tan{\frac{i}{2}} \sin{\Omega} \\ 
    q_2 & = \tan{\frac{i}{2}} \cos{\Omega}  \ .
\end{align}
These elements define the axes of the equinoctial reference frame, $\Sigma_{eq}$, 
\begin{align}
    \Sigma_{eq} = \{O; \ \vect{e}_X, \ \vect{e}_Y, \ \vect{e}_Z \} \ ,
\end{align}
where
\begin{equation}
    \begin{aligned}\label{eq:equinocframe}
        \vect{e}_X & = \frac{1}{1 + q_1^2 + q_2^2} \left[ 1 - q_1^2 + q_2^2, \ \ 2 q_1 q_2, \ \ -2q_1  \right]^T \\
        \vect{e}_Y & = \frac{1}{1 + q_1^2 + q_2^2} \left[ 2 q_1 q_2, \ \ 1 + q_1^2 - q_2^2,  \ \ 2q_2  \right]^T \\
        \vect{e}_Z & = \frac{1}{1 + q_1^2 + q_2^2} \left[ 2 q_1, \ \ -2q_2, \ \ 1 - q_1^2 - q_2^2 \right]^T \ .
    \end{aligned}    
\end{equation}
With the M-GEqOE set defined, the time derivatives for elements are detailed, providing a general form that may accommodate any perturbing forces. 

\subsection{M-GEqOE Time Derivatives}
To obtain the derivatives of the M-GEqOE set with respect to time, projections of the total and external perturbing forces into the orbital reference frame, $\Sigma_{or}$, are necessary. From \cref{eq:orbrefframe} and \cref{eq:pertforces}, components of these projections are given as
\begin{align}
    F_r & = \vect{F} \cdot \vect{e}_r \ , \quad F_f = \vect{F} \cdot \vect{e}_f \ , \quad \text{and} \quad F_h = \vect{F} \cdot \vect{e}_h \ , \\
    P_r & = \vect{P} \cdot \vect{e}_r \ , \quad P_f = \vect{P} \cdot \vect{e}_f \ , \quad \text{and} \quad P_h = \vect{P} \cdot \vect{e}_h \ .
\end{align}
Similarly, the angular velocity of the equinoctial reference frame, $\Sigma_{eq}$, with respect to the inertial frame, $\Sigma$, is projected onto the equinoctial axes as
\begin{align}
    w_X & = F_h \frac{r}{h} \cos{L} \\
    w_Y & = F_h \frac{r}{h} \sin{L} \\
    w_Z & = -F_h \frac{r}{h} \tan{\frac{i}{2}}\sin\left({\omega + \theta}\right) \ .
\end{align}
Other quantities necessary to compute the time derivatives of the M-GEqOEs include the rate of change of the total energy of the system, which is evaluated as
\begin{align}
    \dot{\mathcal{E}} = \frac{\partial U}{\partial t} + \dot{r}P_r + \frac{h}{r} P_f \ .
\end{align}
Then, following the procedure for the GEqOEs detailed in \cite{Bau2021}, the general form of the time derivatives of the M-GEqOEs is obtained as
\begin{equation}
\begin{aligned}\label{eq:geqoeeoms}
    \dot{\tilde{p}} & = \frac{2 \tilde{h}}{\mu_{_{C}}} \left[ 
 \frac{r^2 \dot{\mathcal{E}}}{\tilde{h}} + \frac{r \dot{r}}{\tilde{h}} (2 U - r F_r) \right] \\
    \dot{p}_1 & = p_2 \left( \frac{h - \tilde{h}}{r^2} - w_h \right) + \frac{1}{\tilde{h}} \left( \frac{X}{\tilde{a}} + 2 p_2 \right)(2 U - rF_r) + \frac{1}{\tilde{h}^2} \left[ Y (r + \tilde{p}) + r^2 p_1 \right] \dot{\mathcal{E}} \\
    \dot{p}_2 & = p_1 \left( w_h - \frac{h - \tilde{h}}{r^2} \right)  - \frac{1}{\tilde{h}} \left( \frac{Y}{\tilde{a}} + 2p_1 \right)(2 U - rF_r) + \frac{1}{\tilde{h}^2} \left[ X (r + \tilde{p}) + r^2 p_2 \right] \dot{\mathcal{E}} \\
    \dot{q}_1 & = \frac{1}{2} w_Y \left( 1 + q_1^2 + q_2^2 \right) \\
    \dot{q}_2 & = \frac{1}{2} w_X \left( 1 + q_1^2 + q_2^2 \right) \\
    \dot{L} & = \frac{h}{r^2} + \frac{r}{h} F_h \tan{\left(\frac{i}{2}\right)} \sin{(\omega + \theta)} \ .
\end{aligned}
\end{equation}
All significant perturbations necessary for accurately modeling the dynamics of the object are, thus, directly incorporated into the generalized elements and their associated time derivatives, independent of whether they are expressed as perturbing potentials or forces.

\subsection{Transformations between Cartesian and M-GEqOE Coordinates}

\subsubsection{Cartesian Coordinates to M-GEqOEs}
Let $\vect{x}_C = [ \vect{r} \ \ \dot{\vect{r}} ]^T$ represent the Cartesian state vector comprised of the position and velocity of the object. The transformation from Cartesian into the M-GEqOE coordinates is summarized below:
\begin{itemize}
    \item To obtain the first element, $\tilde{p}$, the generalized angular momentum, $\tilde{h}$, is determined via \cref{eq:genangmom}, and substituted directly into \cref{eq:gensemilatrec}.
    \item To determine the second and third elements, $p_1$ and $p_2$, the generalized eccentricity vector is first computed via \cref{eq:geneccvec}, supplying $p_1$ and $p_2$ as
    \begin{align}
        p_1 = \vect{\tilde{e}} \cdot \vect{e}_Y \quad \text{and} \quad p_2 = \vect{\tilde{e}} \cdot \vect{e}_X \ .
    \end{align}
    \item The fourth and fifth generalized elements, $q_1$ and $q_2$, are determined as
    \begin{align}
        q_1 = \frac{\vect{e}_h \cdot \vect{e}_x}{1 + \vect{e}_h \cdot \vect{e}_z} \quad \text{and} \quad q_2 = - \frac{\vect{e}_h \cdot \vect{e}_y}{1 + \vect{e}_h \cdot \vect{e}_z} \ .
    \end{align}
    \item To obtain the final element, $L$, the position vector in the equinoctial reference frame is found
    \begin{align}
        X = \vect{r} \cdot \vect{e}_X \quad \text{and} \quad Y = \vect{r} \cdot \vect{e}_Y \ .
    \end{align}
    The true longitude may then be determined as
    \begin{align}
        L = \mathrm{atan2}(Y,X) \ ,
    \end{align}
    where $\mathrm{atan2}(\cdot , \cdot)$ denotes the four-quadrant inverse tangent.
\end{itemize}
This transformation from Cartesian coordinates to the M-GEqOE set is utilized to initialize the motion of the object, enabling direct propagation in M-GEqOE coordinates. 

\subsubsection{M-GEqOEs to Cartesian Coordinates}
Let $\vect{x}_M = [\tilde{p} \ \ p_1 \ \ p_2 \ \ q_1 \ \ q_2 \ \ L]^T$ represent the M-GEqOE set at a given time. The transformation from the M-GEqOEs into Cartesian coordinates is detailed below:
\begin{itemize}
    \item The magnitude of the generalized angular momentum vector is determined directly via the first element of the M-GEqOE set as
    \begin{align}
        \tilde{h} = \sqrt{\tilde{p} \mu_{_{C}}} \ .
    \end{align}
    \item Next, the elements $p_1, \ p_2,$ and $L$ yield the magnitude of the position vector in the equinoctial frame as
    \begin{align}
        r = \frac{\tilde{p}}{1 + p_1 \sin{L} + p_2 \cos{L}} \ .
    \end{align}
    \item The components of the position vector in the equinoctial frame are determined via
    \begin{equation}
        \begin{aligned}
            X & = r \cos{L} \\
            Y & = r \sin{L} \ .
        \end{aligned} 
    \end{equation}
    \item Leveraging \cref{eq:equinocframe}, the elements $q_1$ and $q_2$ are used to determine the unit vectors that define the equinoctial frame, $[\vect{e}_X, \ \vect{e}_Y, \ \vect{e}_Z]^T$. 
    \item The position vector is then found as
    \begin{align}
        \vect{r} = X \vect{e}_X + Y \vect{e}_Y \ .
    \end{align}
    \item Next, the radial velocity, $\dot{r}$, is determined as
    \begin{align}
        \dot{r} = \frac{\mu_{_{C}}}{\tilde{h}} \left( p_2 \sin{L} - p_1 \cos{L} \right) \ .
    \end{align}
    \item Finally, the velocity vector is expressed as
    \begin{align}
        \dot{\vect{r}} = \dot{X} \vect{e}_X + \dot{Y} \vect{e}_Y \ ,
    \end{align}
    where
    \begin{equation}
        \begin{aligned}
            \dot{X} & = \frac{\dot{r} X}{r} - \frac{h Y}{r^2} \\
            \dot{Y} & = \frac{\dot{r} Y}{r} + \frac{h X}{r^2} \ .
        \end{aligned}
    \end{equation}
\end{itemize}
The transformation from M-GEqOEs into Cartesian coordinates is utilized to visualize spacecraft trajectories in various frames and to assess the accuracy of propagating the M-GEqOE equations of motion against conventional methods.

\section{Cislunar Dynamics Modeling}
As spacecraft and other objects traverse cislunar space, their complex dynamical motion is influenced by the gravitational forces of both the Earth and the Moon. Therefore, in addition to the dominant gravitational influence of one body, the gravitational perturbation caused by the other must be accounted for in trajectory design and prediction. In addition, other perturbations, such as the gravitational influence of the Sun and the oblateness effects of the Earth and Moon, may also be significant. In the current work, a high-fidelity simulation of cislunar dynamics is achieved by modeling the gravitational influence of the Earth, Moon, and Sun on the object. The formulation of the dynamics is detailed below in the context of both Cartesian and M-GEqOE coordinates. While the discussion pertains only to the three celestial bodies, the framework may be extended to incorporate additional point-mass gravitational effects or other perturbations.

\subsection{Cartesian Coordinates}\label{sec:cartdyn}
The influence of the Earth, Moon, and Sun is modeled via an $N$-body ephemeris model that is based in an inertial frame centered on the primary gravitational body. Assuming that the mass of the object of interest is negligible relative to the masses of the other bodies, the general form for the equations of motion that govern the object dynamics is expressed as the second-order relative vector differential equations 
\begin{align}
    \Ddot{\vect{r}}_{{C s/c}} = -\frac{GM_{C}}{r_{{C s/c}}^3}\vect{r}_{{C i}} + G \sum_{\substack{i=1\\
                  i\neq C,s/c}}^{N} m_i \left(\frac{\vect{r}_{{s/c i}}}{r_{{s/c i}}^3} - \frac{\vect{r}_{{C i}}}{r_{{C i}}^3} \right) \ ,
\end{align}
where $\vect{r}_{{C s/c}}$ represents the position of the object relative to the central body, $\vect{r}_{s/c i}$ denotes the position of each perturbing body with respect to the object, and $\vect{r}_{C i}$ locates the position of each perturbing body relative to the central body. The DE440 planetary ephemerides retrieved from the NASA Jet Propulsion Lab (JPL) Navigation and Ancillary Information Facility (NAIF) are employed to obtain the relevant quantities in this model \citep{DE440}. Finally, to mitigate the accumulation of numerical error during propagation, the equations of motion are nondimensionalized. One convenient nondimensionalization scheme for this system employs characteristic quantities associated with the Earth–Moon Circular Restricted Three-Body Problem (CR3BP) \citep{Broucke1968}. The characteristic length, $l^*$, is chosen as the mean Earth-Moon distance. The characteristic mass, $m^*$, is defined as the sum of the masses of the Earth and the Moon. The characteristic time, $t^*$, is then defined as
\begin{align}
    t^* = \sqrt{\frac{{l^*}^3}{G m^*}} \ .
\end{align}
Together, these characteristic quantities are employed in the nondimensional propagation of the motion of the object in Cartesian coordinates. 

\subsection{M-GEqOE Coordinates}\label{sec:geqdyn}
For high-fidelity cislunar propagation using the M-GEqOEs, the gravitational influences of the Earth, Moon and Sun are included. When the Earth is treated as the central gravitational body, the Moon is modeled as a point-mass third-body perturbation. Conversely, for trajectories in the vicinity of the Moon, the Moon is treated as the central body, with the Earth providing the perturbing potential. For the case of a single additional perturbing body, the general form of the gravitational potential energy is given by
\begin{align}\label{eq:3bpert}
    U_{{P}} = \mu_{P} \left( \frac{1}{r_{{P sc}}} - \frac{\vect{r} \cdot \vect{r}_{{CP}}}{r_{{CP}}^3} \right) \ ,
\end{align}
where $\vect{r}$ denotes the position vector from the central body to the spacecraft, and the subscripts $C$ and $P$ refer to the central and perturbing bodies, respectively \citep{Battin1999}. Although the perturbing potential may be approximated using a convergent series of Legendre polynomials \citep{GuptaDeMars2025Astro}, the exact form in \cref{eq:3bpert} exhibits closer agreement with the Cartesian formulation. In addition to third-body perturbing effects, the point-mass gravitational effects of the Sun are also considered in the current work. However, noting the numerical instabilities resulting from the use of true solar ephemerides, solar gravity is incorporated as an external perturbing force, denoted by $\vect{P}$ in \cref{eq:pertforces} \citep{Bau2021}. Nevertheless, these effects do manifest directly in the equations that evolve the M-GEqOEs, as indicated in \cref{eq:geqoeeoms}. Finally, consistent with the Cartesian formulation, the characteristic quantities of the Earth–Moon CR3BP are used to nondimensionalize the M-GEqOE propagation, thereby mitigating the accumulation of numerical error.

\subsubsection{Constant Potential Offset}
In the M-GEqOE formulation employed in the current work, the quantities associated with generalized orbital motion depend on the effective potential, defined in \cref{eq:Ueff}. Of particular note is the generalized angular momentum, $\tilde{{h}}$ defined in \cref{eq:genangmom}, which must remain nonnegative by definition. When the third-body gravitational potential is evaluated directly, its nominal zero point can result in a negative $U_{\mathrm{eff}}$, even though the underlying dynamics remain physically valid. In such cases, the numerical evaluation of $\tilde{h}$ becomes ill-defined, which in turn leads to failures in the computation of the M-GEqOEs.

As such, to ensure the validity of $U_{\mathrm{eff}}$ throughout the propagation, a constant offset $U_{\mathrm{offset}}$ is applied to the third-body gravitational potential. An instantaneous offset is first evaluated along the Cartesian trajectory by enforcing the nonnegativity of $U_\mathrm{eff}$ as
\begin{align}
    U_{\mathrm{offset}}(t) = -U(t) - \frac{h(t)^2}{2 r(t)^2} - \frac{\dot{r}(t)^2}{4} + \frac{\mu_{_{C}}}{2r(t)} \ .
\end{align}
To ensure that $U_{\mathrm{eff}}$ remains nonnegative for all times, the constant offset used for the propagation is then selected as
\begin{align}
    U_{\mathrm{offset}} = \max_{t \in T} U_{\mathrm{offset}}(t) \ ,
\end{align}
where $T$ represents the total propagation time. This offset corresponds to a shift in the reference level of the potential energy and does not alter the equations of motion, since the force model depends only on the gradient of the potential. 

\section{Uncertainty Characterization}
For a chaotic regime such as the cislunar dynamical environment, where small perturbations can manifest into vastly different trajectories, realistic propagation of uncertainty is crucial for SDA. In the current work, the propagation and characterization of uncertainty in cislunar space are examined in both M-GEqOE and Cartesian coordinate representations to evaluate the influence of coordinate choice on uncertainty realism. Monte Carlo simulations are employed to characterize the growth of uncertainty along various representative orbits. In previous work, statistical measures such as the Kullback-Leibler divergence \citep{GuptaDeMars2025ESA,GuptaDeMars2025Unc} and the Bhattacharyya coefficient \citep{GuptaDeMars2025AAS} are then employed in order to characterize downstream uncertainty. However, these measures require explicit specification of a true and an approximating distribution, of which the latter is assumed to be Gaussian. In the case where the true distribution is unknown and only accessible through samples, these measures introduce assumptions that obscure whether deviations arise from intrinsic non-Gaussian structure or from the choice of approximating distribution. As such, a different measure of multivariate normality, the Henze--Zirkler (HZ) test, is explored \citep{Henze1990}. The HZ test operates directly on the Monte Carlo sample statistics and their pairwise distances, without any dependence on an assumed true pdf or an approximating pdf. The HZ test is applied to identify the onset of non-Gaussianity for uncertainty propagated in each coordinate set. 

\subsection{Monte Carlo Analysis}
In the current work, Monte Carlo analysis utilizes $N = 10,000$ samples for all simulations to facilitate analysis while maintaining computational tractability. Uncertainty propagation in the M-GEqOE framework follows the high-fidelity dynamics described in \cref{sec:geqdyn}, including third-body perturbations and solar gravity effects. Cartesian samples are propagated using the full $N-$body ephemeris equations. In both formulations, the ensemble states are recorded at each time step to assess the characteristics of the underlying probability distribution. 

\subsection{Henze--Zirkler Test for Multivariate Normality}
In order to identify the emergence of non-Gaussianity in uncertainty propagation, the Henze--Zirkler (HZ) test is employed. The HZ test falls under the category of \textit{consistent} approaches, implying that the test should consistently reject all non-multivariate normal distributions \citep{Henze1990}. While other tests for multivariate normality exist, the HZ test has been demonstrated to provide a balance between statistical power and computational efficiency \citep{Mecklin2005}. 

To implement the HZ test, let $\vect{x}_i (t) \in \mathbb{R}^n$ for $i = 1, 2, ... , N$ denote an ensemble of particles that represent the probability distribution at time $t$. The sample mean, $\vect{m}_{x}(t)$, and covariance, $\vect{P}_{xx}(t)$, are computed as
\begin{align}
    \vect{m}_x (t) & = \frac{1}{N} \sum_{i=1}^{N} \vect{x}_i (t) \\
    \vect{P}_{xx} (t) & = \frac{1}{N} \sum_{i=1}^{N} \left (\vect{x}_i (t) - \vect{m}_x(t) \right) \left (\vect{x}_i (t) - \vect{m}_x(t) \right)^T \ .
\end{align}
Next, the squared Mahalanobis distances between each sample and the mean, as well as between two samples, are determined as
\begin{align}
    d_i^2 (t) & = \left( \vect{x}_i(t) - \vect{m}_x(t) \right)^T \vect{P}_{xx}^{-1} \left( \vect{x}_i(t) - \vect{m}_x(t) \right) \\
    d_{ij}^2 (t) & = \left( \vect{x}_i(t) - \vect{x}_j(t) \right)^T \vect{P}_{xx}^{-1} \left( \vect{x}_i(t) - \vect{x}_j(t) \right) \ .
\end{align}
A smoothing parameter, denoted $\beta$, that depends on the sample size, $N$, and the state dimension, $n$, is selected such that
\begin{align}
    \beta = \frac{1}{\sqrt{2}} \left( \frac{N (2n + 1)}{4} \right)^{\frac{1}{n+4}} \ .
\end{align}
The HZ statistic is then determined as
\begin{align}
    \begin{split}
        \mathrm{HZ}(t) = {}& \left[ \frac{1}{N^2} \sum_{i=1}^{N} \sum_{j=1}^{N} \exp{ \left\{ - \frac{\beta^2 d_{ij}^2(t)}{2}  \right\} } \right] - \left[ \frac{2 \gamma^{-n/2}}{N} \sum_{i=1}^N \exp{ \left\{ - \frac{\beta^2 d_i^2(t)}{2\gamma} \right\} } \right] \\
    & + \left[ (1 + 2 \beta^2)^{-n/2} \right] \ ,
    \end{split}
\end{align}
where $\gamma = 1 + \beta^2$. The null hypothesis, $H_0$, states that the samples are drawn from a multivariate normal distribution. Under $H_0$, the distribution of the HZ statistic is approximated by a log-normal distribution \citep{Flegel2017}. In the current work, hypothesis testing is conducted following the methodology of \cite{Ortiz2026}. The $p$-value is obtained from the cumulative distribution function of the corresponding log-normal distribution as
\begin{align}
    p(t) = 1 - \mathrm{lognormal}\left( \mathrm{HZ}(t), \hat{\mu}(t), \hat{\sigma}(t) \right) \ ,
\end{align}
where $\hat{\mu}(t)$ and $\hat{\sigma}(t)$ denote the mean and standard deviation of the log-normal distribution associated with the HZ statistic \citep{Flegel2017}. Hypothesis rejection is then determined relative to a given significance level, denoted $\alpha$, as
\begin{align}
    p & > \alpha \ \rightarrow \ \ \ H_0 \ \mathrm{cannot \ be \ rejected} \\
    p & \leq \alpha \ \rightarrow \ \ \ H_0 \ \mathrm{should \ be \ rejected} \ .
\end{align}
In the current work, a significance level of $\alpha = 0.05$ is selected. Together, the HZ statistic and its associated $p$-value are evaluated along each trajectory to assess how the validity of a Gaussian uncertainty assumption evolves in time, and whether the choice of the coordinate representation used for uncertainty propagation influences the emergence of non-Gaussian behavior.

\section{Results}
State and uncertainty propagation under high-fidelity dynamics are assessed for a variety of cislunar orbits that traverse various regions within this domain. For each orbit, the states are evolved in generalized coordinates and, via transformation into Cartesian coordinates, compared against the solutions obtained via direct Cartesian propagation. In addition to state propagation, the evolution of uncertainty along each trajectory is assessed for both sets of coordinates. To maintain the scope of this work as the characterization of downstream probability distributions, a single revolution of each cislunar trajectory is considered for uncertainty propagation.

\subsection{9:2 Near Rectilinear Halo Orbit}
\begin{figure}[!htbp]
    \centering
    \subfloat[Earth-Moon rotating frame]
        {\includegraphics[width=0.48\textwidth]{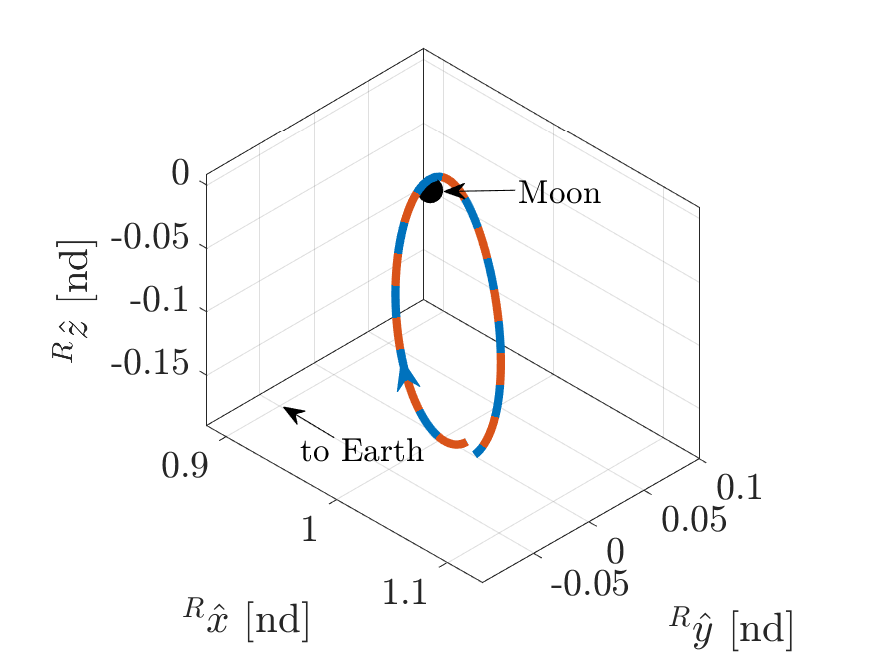}
        \label{fig:nrhorot}}
    \subfloat[Moon-centered inertial frame]
        {\includegraphics[width=0.48\textwidth]{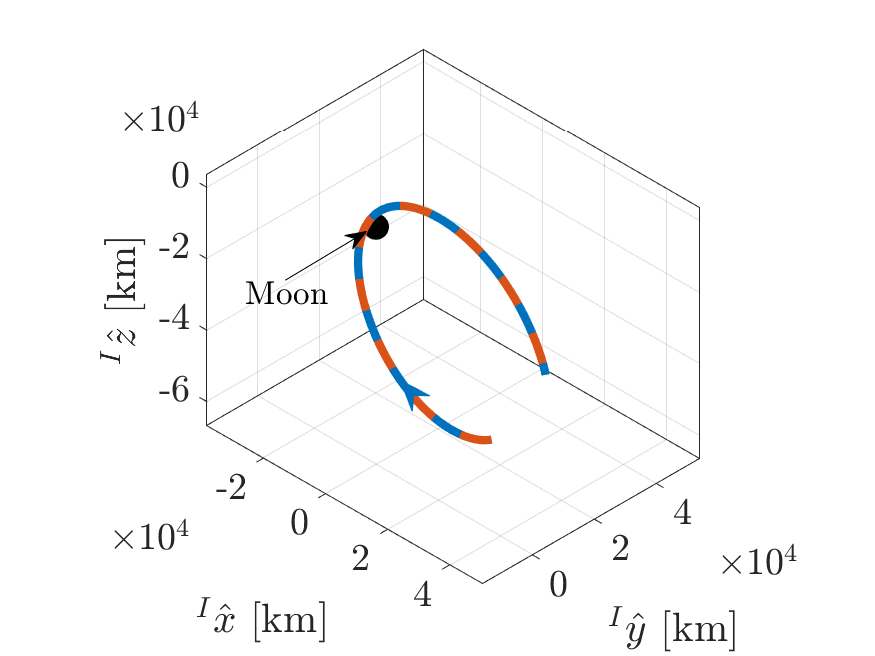}
        \label{fig:nrhomci}}   \\
    \subfloat[Position and velocity errors]
        {\includegraphics[width=0.6\textwidth]{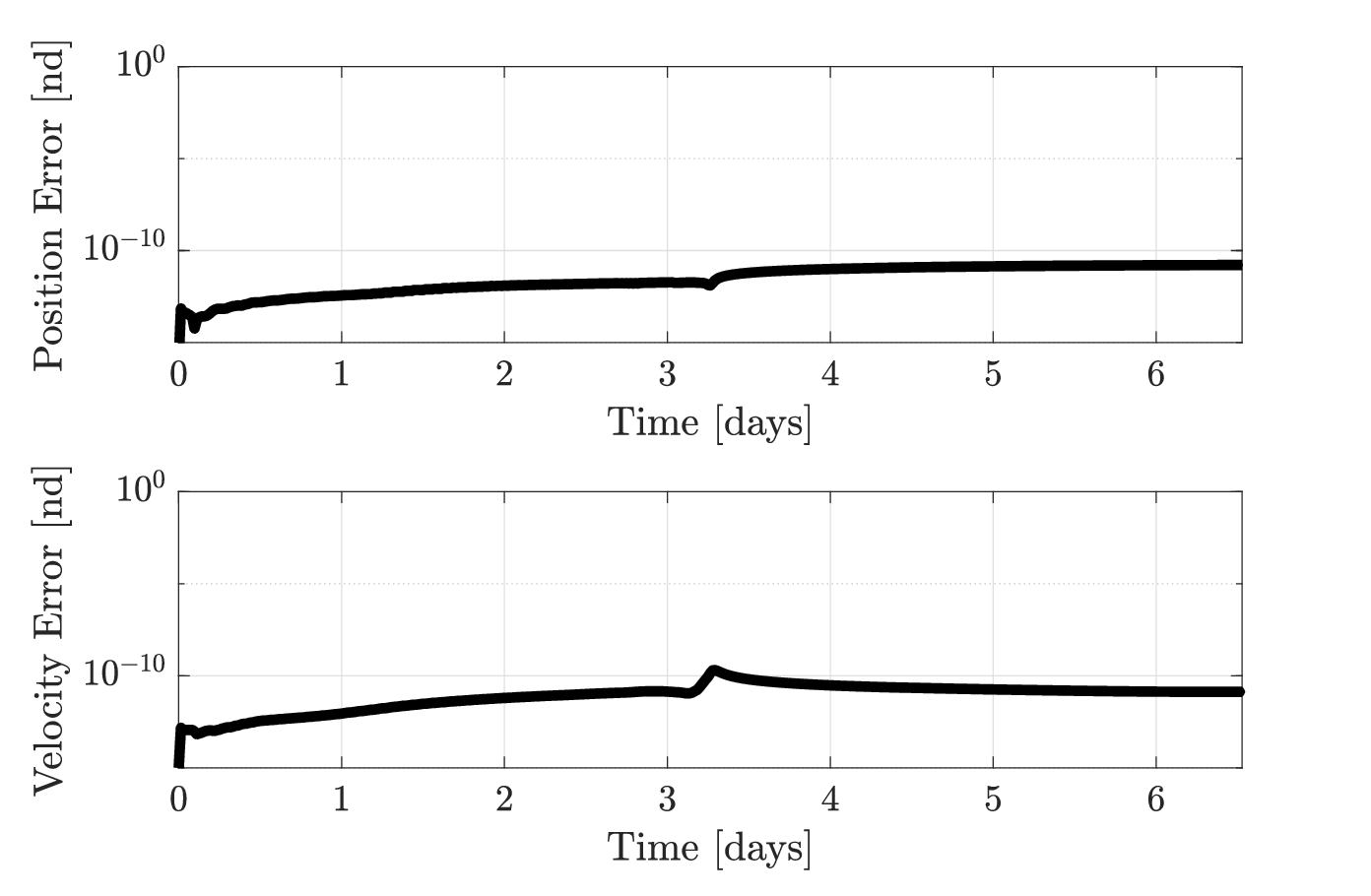}
        \label{fig:nrhoerror}}   \\    
    \subfloat[Evolution of M-GEqOEs]
    {\includegraphics[width=0.9\textwidth]{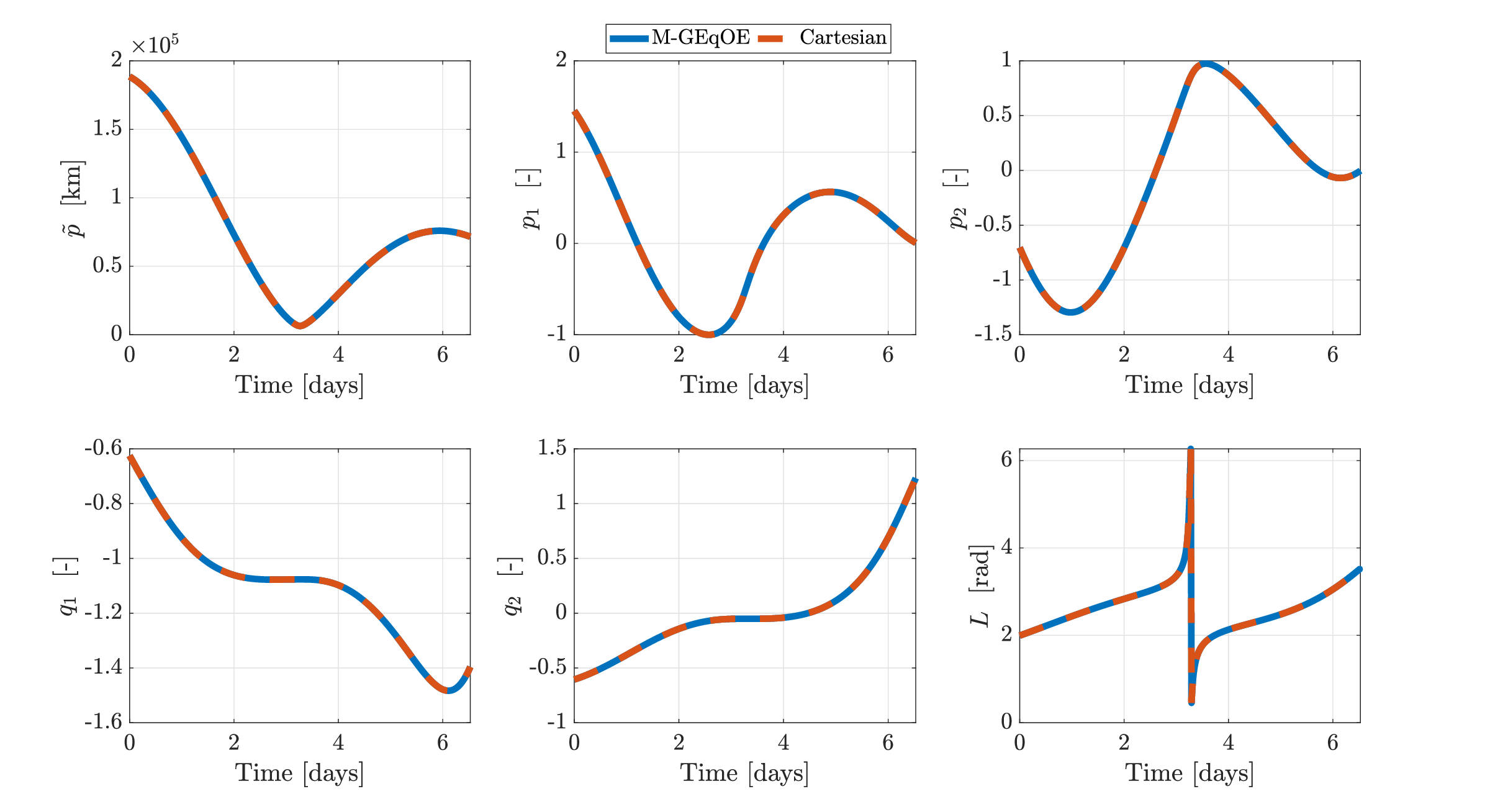}\label{fig:nrhomgeqoe}}
    \caption{Cartesian and M-GEqOE representations of the 9:2 NRHO over $6.5 \ days$. Blue and orange curves represent the M-GEqOE and Cartesian solutions, respectively.}
    \label{fig:nrhostate}
\end{figure}

The \rt{9}{2} Near Rectilinear Halo Orbit (NRHO), which is the baseline orbit for NASA's upcoming Gateway lunar outpost \citep{Lee2019}, is analyzed. The dynamics that govern the motion of a spacecraft in this orbit are highly nonlinear, and unmodeled perturbations may induce significant deviation from the baseline trajectory over time. Note that, in this case, the resonance ratio indicates synodic resonance with the Moon, a characteristic that allows eclipse feasibility for both Gateway and the Orion spacecraft \citep{Zimovan-Spreen2022}. The NRHO is propagated using the M-GEqOE equations for one revolution of the orbit or approximately $6.5 \ days$. The resulting M-GEqOE solution is mapped into Cartesian coordinates, as plotted in blue in the Earth-Moon rotating and Moon-centered inertial frames in \cref{fig:nrhorot} and \cref{fig:nrhomci}, respectively. The trajectory obtained via direct propagation in Cartesian coordinates also appears in these frames in orange. While the Cartesian representations of the two solutions indicate visual consistency between the results, the error in position and velocity confirms the accuracy of the M-GEqOE solution, plotted in \cref{fig:nrhoerror}. The error in position and velocity between the two methods of propagation remains generally bounded, with an increase of approximately one order of magnitude observed at the perilune crossing, about $3.25 \ days$ into propagation. Finally, mapping the Cartesian solution into the M-GEqOE coordinates confirms the accuracy of the transformation over time, as illustrated in \cref{fig:nrhomgeqoe}. 

\begin{figure}[!htbp]
    \centering
    \subfloat[HZ statistic]
        {\includegraphics[width=0.4925\textwidth]{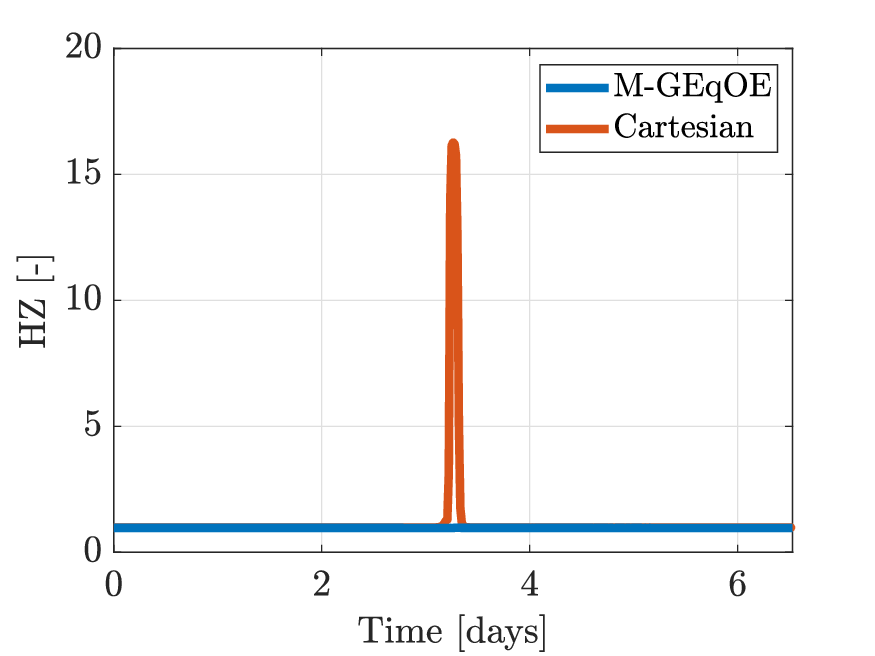}
        \label{fig:nrhohz}}
    \subfloat[p-value]
        {\includegraphics[width=0.4925\textwidth]{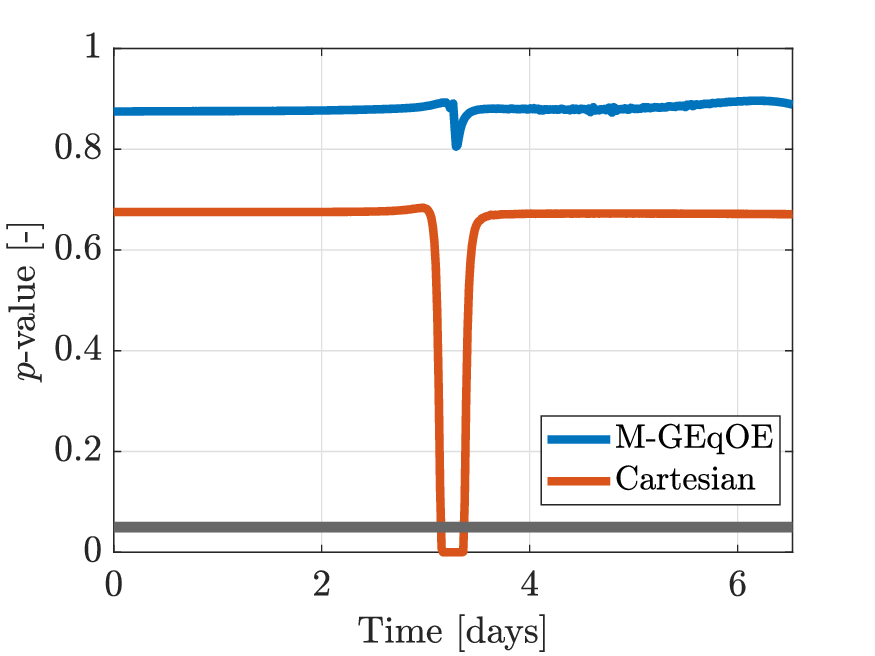}
        \label{fig:nrhopval}}  
    \caption{Henze--Zirkler test applied to uncertainty propagated along the 9:2 NRHO.}
    \label{fig:nrhounc}
\end{figure}

Next, the evolution of uncertainty along the NRHO is assessed for the M-GEqOEs and compared against propagation in Cartesian coordinates. At the initial time, which occurs at apolune, $1\sigma$ values of uncertainty of $1 \ km$ and $1 \ cm/s$ are assumed along the position and velocity channels, respectively. Uncertainty is then evaluated at each point along the orbit for approximately $6.5 \ days$. \cref{fig:nrhounc} illustrates the HZ test statistic and the associated $p$-value evaluated along the NRHO for both the M-GEqOE (blue) and Cartesian (orange) representations. Values of the HZ statistic exceeding unity indicate increasing departure from Gaussian behavior, while $p$-values below the $0.05$ threshold, denoted by the grey bound in \cref{fig:nrhopval}, correspond to rejection of the Gaussian hypothesis. The Cartesian methodology exhibits a sharp increase in the HZ statistic, which is accompanied by a pronounced decrease in the $p$-value, at the time of perilune passage along the orbit. For the M-GEqOE methodology, consistently low values of the HZ statistic are observed, with $p$-values that remain above the rejection threshold, for the full trajectory. As such, uncertainty evolved in the M-GEqOE coordinates better preserves Gaussian behavior over time.

\begin{figure}[!htbp]
    \centering
    {
    \begin{tikzpicture}
        \node (popfig) at ([xshift=0cm,yshift=0cm] current page.center) {\includegraphics[width=1.05\textwidth]{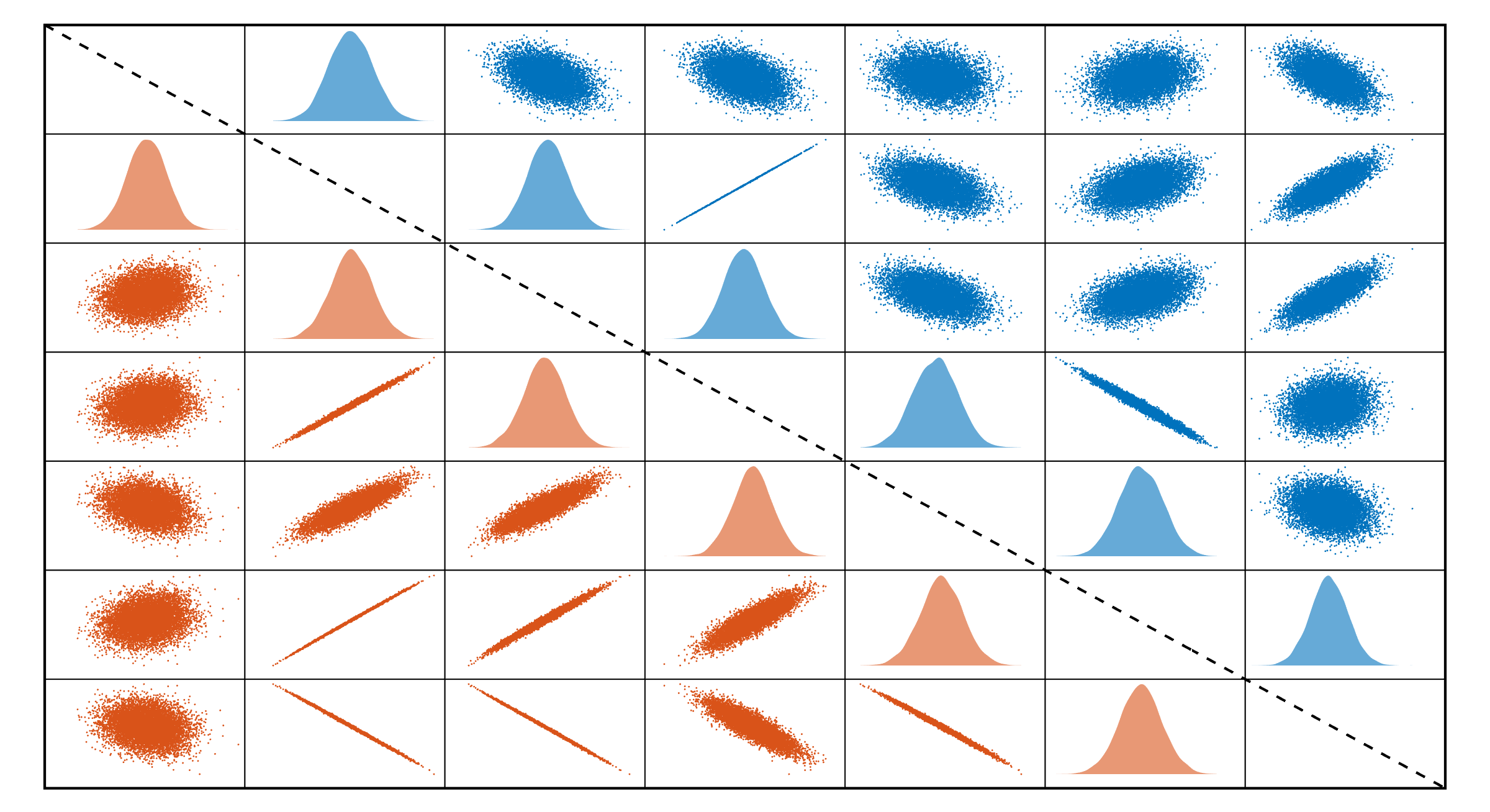}};
        \node at ([xshift=-6.75cm,yshift=1.0cm] popfig.center) {\large {$y$}};
        \node at ([xshift=-6.75cm,yshift=0cm] popfig.center) {\large {$z$}};
        \node at ([xshift=-6.75cm,yshift=-1.0cm] popfig.center) {\large {$\dot{x}$}};
        \node at ([xshift=-6.75cm,yshift=-2.0cm] popfig.center) {\large {$\dot{y}$}};
        \node at ([xshift=-6.75cm,yshift=-3.0cm] popfig.center) {\large {$\dot{z}$}};
        \node at ([xshift=-5.4cm,yshift=-3.9cm] popfig.center) {\large {$x$}};
        \node at ([xshift=-3.6cm,yshift=-3.9cm] popfig.center) {\large {$y$}};
        \node at ([xshift=-1.8cm,yshift=-3.9cm] popfig.center) {\large {$z$}};
        \node at ([xshift=0cm,yshift=-3.9cm] popfig.center){\large {$\dot{x}$}};
        \node at ([xshift=1.8cm,yshift=-3.9cm] popfig.center) {\large {$\dot{y}$}};
        \node at ([xshift=3.6cm,yshift=-3.9cm] popfig.center) {\large {$\dot{z}$}};
        \node at ([xshift=-3.6cm,yshift=3.9cm] popfig.center) {\large {$\tilde{p}$}};
        \node at ([xshift=-1.8cm,yshift=3.9cm] popfig.center) {\large {$p_1$}};
        \node at ([xshift=0cm,yshift=3.9cm] popfig.center) {\large {$p_2$}};
        \node at ([xshift=1.8cm,yshift=3.9cm] popfig.center) {\large {$q_1$}};
        \node at ([xshift=3.6cm,yshift=3.9cm] popfig.center) {\large {$q_2$}};
        \node at ([xshift=5.4cm,yshift=3.9cm] popfig.center) {\large {$L$}};
        \node at ([xshift=6.8cm,yshift=3cm] popfig.center) {\large {$\tilde{p}$}};
        \node at ([xshift=6.8cm,yshift=2cm] popfig.center) {\large {$p_1$}};
        \node at ([xshift=6.8cm,yshift=1cm] popfig.center) {\large {$p_2$}};
        \node at ([xshift=6.8cm,yshift=0cm] popfig.center) {\large {$q_1$}};
        \node at ([xshift=6.8cm,yshift=-1cm] popfig.center) {\large {$q_2$}};
    \end{tikzpicture}}
    \caption{Pairs plot for the NRHO at the first perilune pass (${t = 3.25 \ days}$) showing projections in Cartesian (lower triangular) and M-GEqOE (upper triangular) coordinates.}
    \label{fig:nrhopop}
\end{figure}

\begin{figure}[!htbp]
    \centering
    {
    \begin{tikzpicture}
        \node (popfig) at ([xshift=0cm,yshift=0cm] current page.center) {\includegraphics[width=1.05\textwidth]{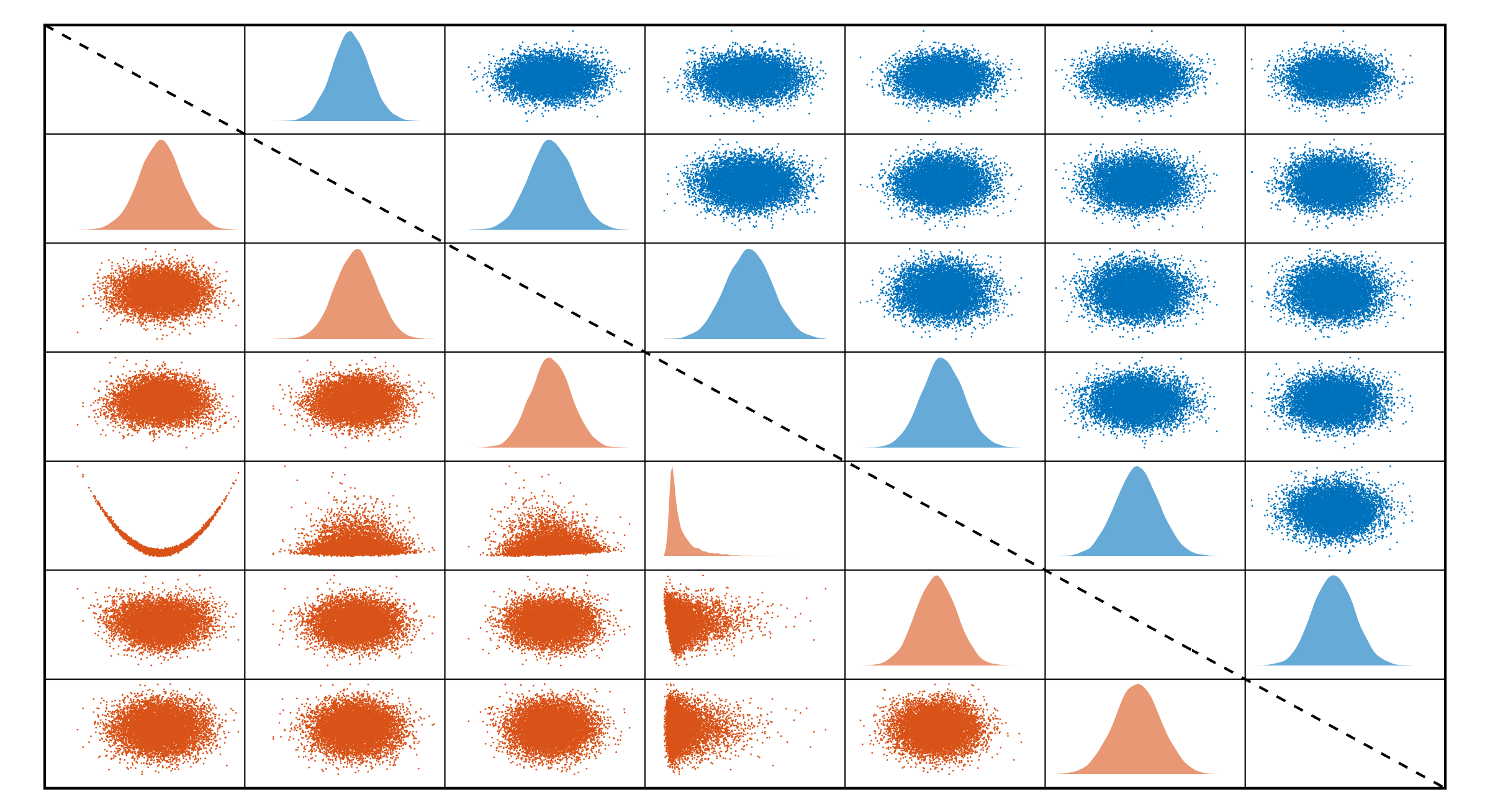}};
        \node at ([xshift=-6.75cm,yshift=1.0cm] popfig.center) {\large {$\lambda_{y}$}};
        \node at ([xshift=-6.75cm,yshift=0cm] popfig.center) {\large {$\lambda_{z}$}};
        \node at ([xshift=-6.75cm,yshift=-1.0cm] popfig.center) {\large {$\lambda_{\dot{x}}$}};
        \node at ([xshift=-6.75cm,yshift=-2.0cm] popfig.center) {\large {$\lambda_{\dot{y}}$}};
        \node at ([xshift=-6.75cm,yshift=-3.0cm] popfig.center) {\large {$\lambda_{\dot{z}}$}};
        \node at ([xshift=-5.4cm,yshift=-3.9cm] popfig.center) {\large {$\lambda_{x}$}};
        \node at ([xshift=-3.6cm,yshift=-3.9cm] popfig.center) {\large {$\lambda_{y}$}};
        \node at ([xshift=-1.8cm,yshift=-3.9cm] popfig.center) {\large {$\lambda_{z}$}};
        \node at ([xshift=0cm,yshift=-3.9cm] popfig.center) {\large {$\lambda_{\dot{x}}$}};
        \node at ([xshift=1.8cm,yshift=-3.9cm] popfig.center) {\large {$\lambda_{\dot{y}}$}};
        \node at ([xshift=3.6cm,yshift=-3.9cm] popfig.center) {\large {$\lambda_{\dot{z}}$}};
        \node at ([xshift=-3.6cm,yshift=3.9cm] popfig.center) {\large {$\lambda_{\tilde{p}}$}};
        \node at ([xshift=-1.8cm,yshift=3.9cm] popfig.center) {\large {$\lambda_{p_1}$}};
        \node at ([xshift=0cm,yshift=3.9cm] popfig.center) {\large {$\lambda_{p_2}$}};
        \node at ([xshift=1.8cm,yshift=3.9cm] popfig.center) {\large {$\lambda_{q_1}$}};
        \node at ([xshift=3.6cm,yshift=3.9cm] popfig.center) {\large {$\lambda_{q_2}$}};
        \node at ([xshift=5.4cm,yshift=3.9cm] popfig.center) {\large {$\lambda_{L}$}};
        \node at ([xshift=6.8cm,yshift=3cm] popfig.center) {\large {$\lambda_{\tilde{p}}$}};
        \node at ([xshift=6.8cm,yshift=2cm] popfig.center) {\large {$\lambda_{p_1}$}};
        \node at ([xshift=6.8cm,yshift=1cm] popfig.center) {\large {$\lambda_{p_2}$}};
        \node at ([xshift=6.8cm,yshift=0cm] popfig.center) {\large {$\lambda_{q_1}$}};
        \node at ([xshift=6.8cm,yshift=-1cm] popfig.center) {\large {$\lambda_{q_2}$}};
    \end{tikzpicture}}
    \caption{Eigenspace pairs plot for the NRHO at the first perilune pass (${t = 3.25 \ days}$) showing projections in Cartesian (lower triangular) and M-GEqOE (upper triangular) coordinates.}
    \label{fig:nrhoeigpop}
\end{figure}
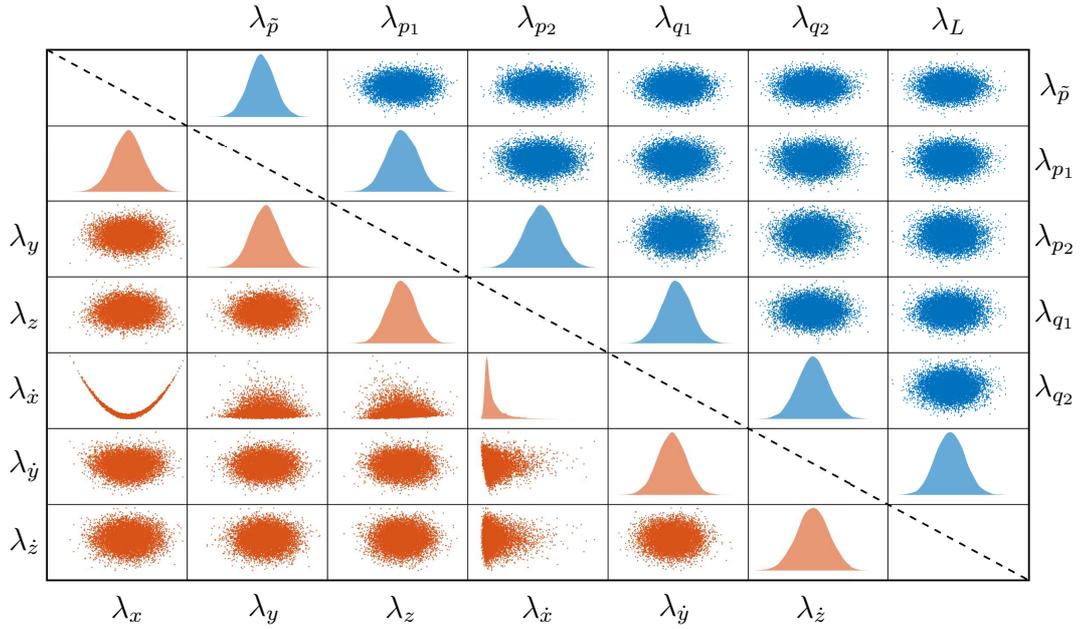

This behavior is further analyzed by visualizing the uncertainty clouds at perilune for both coordinates. \cref{fig:nrhopop} illustrates pairs plots for the M-GEqOE (blue, upper triangle) and Cartesian (orange, lower triangle) methodologies. At this time, the value of the HZ statistic indicates a significant departure from Gaussianity. However, this behavior is not immediately evident from the pairwise projections plotted in \cref{fig:nrhopop}. As such, in order to assess the structure of the particle cloud, the ensemble is transformed into the eigenspace of the covariance matrix. The resulting pairwise projections appear in \cref{fig:nrhoeigpop} for both M-GEqOE (blue) and Cartesian (orange) coordinates. In this view, the departure from Gaussianity is evident, characterized by the curvature observed in the Cartesian projections and reflecting the effects of nonlinearity at this location in orbit. The M-GEqOE representation, however, more closely resembles Gaussian behavior through the perilune passage, without the emergence of any curvature or tail-like structures in the projections.

\subsection{4:1 Sidereal Resonant Orbit}
Orbits with periods commensurate with the lunar sidereal period offer a variety of options for trajectories that traverse the cislunar domain extensively. These orbits are characterized by their resonance ratio, denoted \rt{p}{q}, where a spacecraft orbits the Earth $p$ times in the time that the Moon completes $q$ orbits about the Earth. A sample sidereal resonant orbit considered in the current work is identified from the family of \rt{4}{1} Earth-Moon resonant orbits and possesses a period of approximately $27.30 \ days$ in the rotating frame. Four revolutions of the orbit as viewed in the inertial frame represent one ``closed'' orbit in the rotating frame. Despite being centered on the Earth, this orbit possesses a relatively high apogee radius and, thus, perturbations from lunar gravity are non-negligible. 

\begin{figure}[!htbp]
    \centering
    \subfloat[Earth-Moon rotating frame]
        {\includegraphics[width=0.48\textwidth]{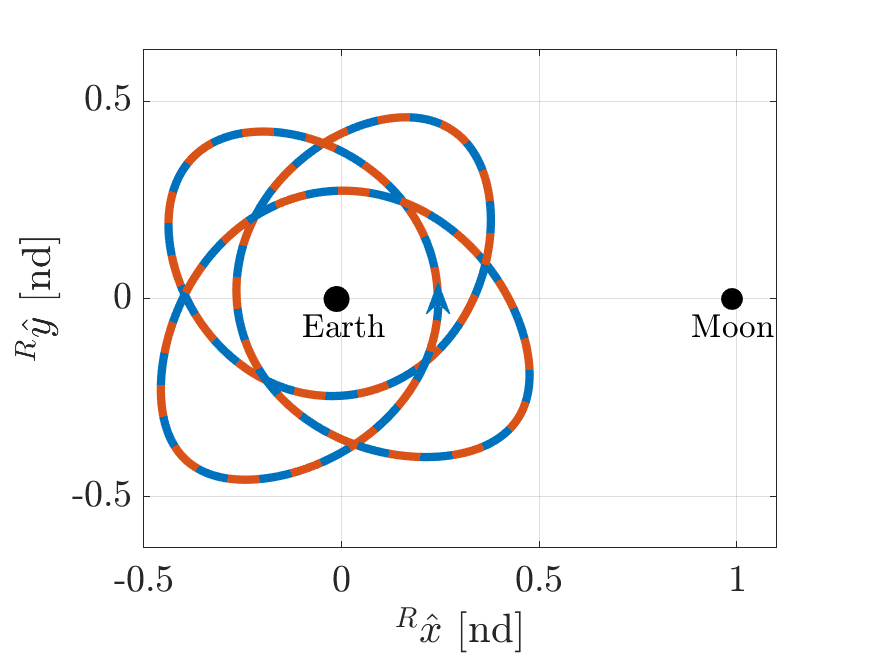}
        \label{fig:r41rot}}
    \subfloat[Earth-centered inertial frame]
        {\includegraphics[width=0.48\textwidth]{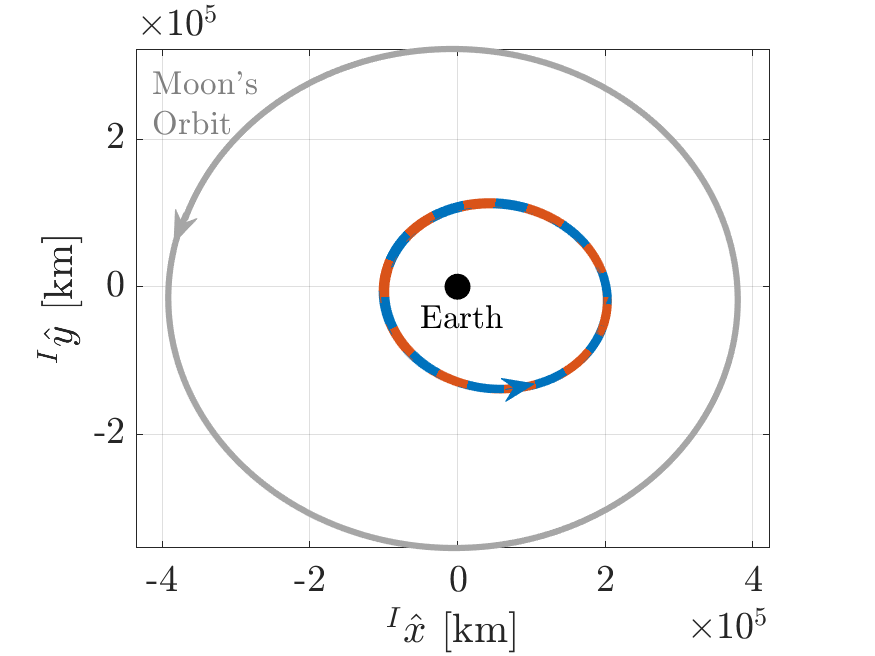}
        \label{fig:r41eci}}   \\
    \subfloat[Position and velocity errors]
        {\includegraphics[width=0.6\textwidth]{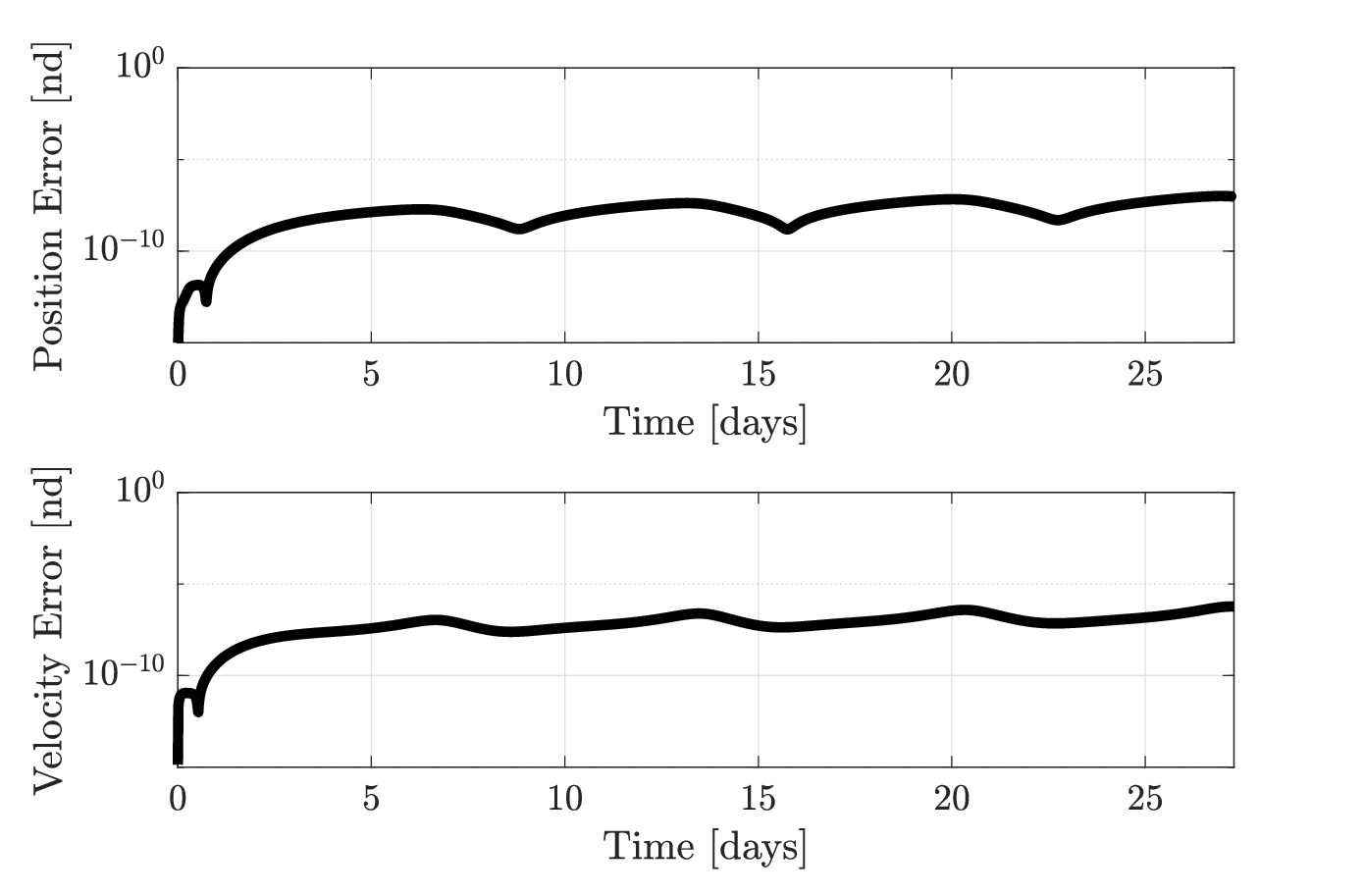}
        \label{fig:r41error}}   \\    
    \subfloat[Evolution of M-GEqOEs]
    {\includegraphics[width=0.9\textwidth]{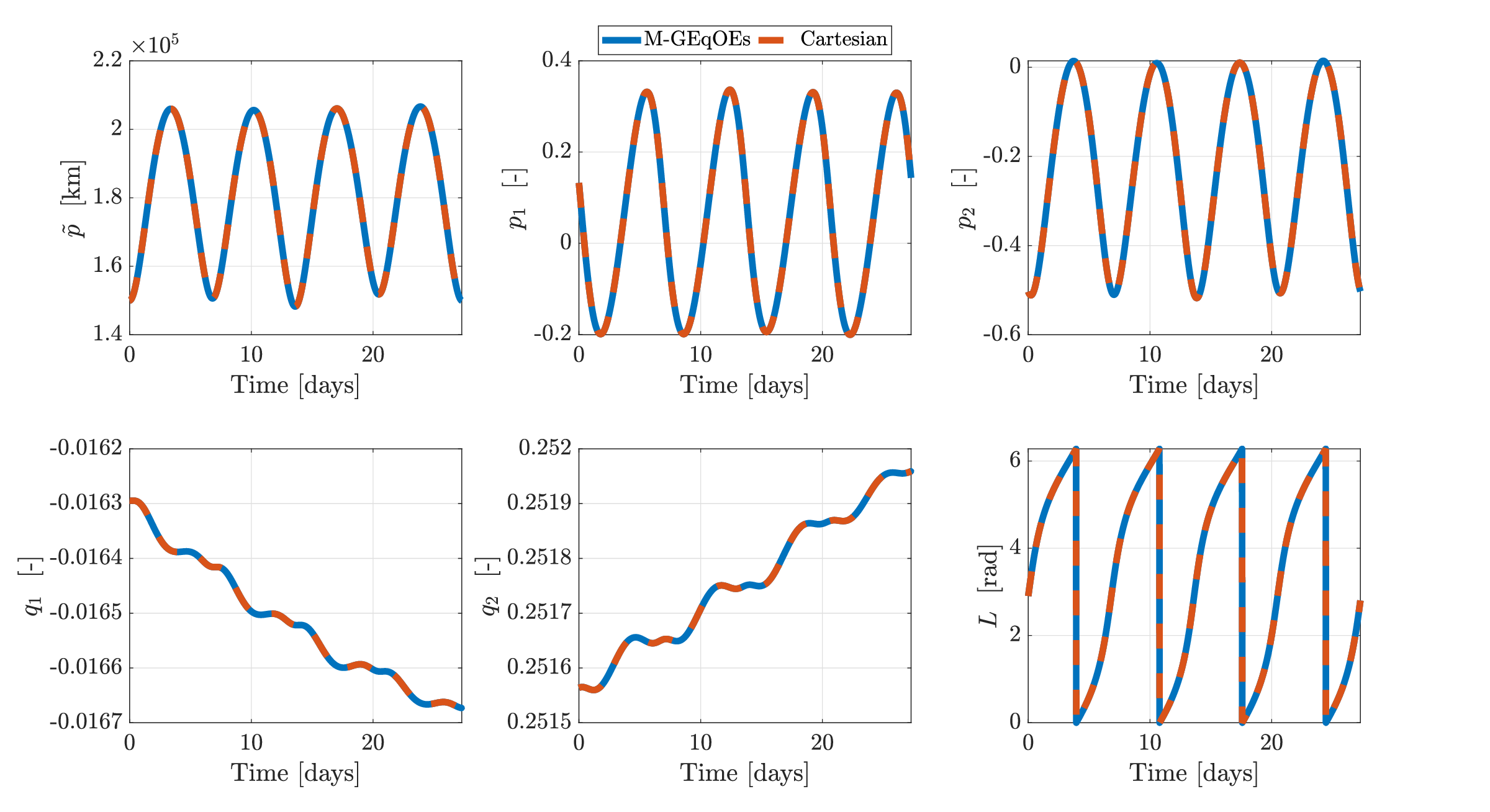}\label{fig:r41mgeqoe}}
    \caption{Cartesian and M-GEqOE representations of the 4:1 resonant orbit over $27.30 \ days$. Blue and orange curves represent the M-GEqOE and Cartesian solutions, respectively.}
    \label{fig:r41state}
\end{figure}

For this analysis, the initial states for the orbit are evolved under Earth-Moon-Sun dynamics using the M-GEqOE and Cartesian equations of motion. The solution obtained via the M-GEqOE methodology is transformed into Cartesian coordinates and appears in \cref{fig:r41rot,fig:r41eci} as viewed in the Earth-Moon rotating and Earth-centered inertial frames. The resulting blue solution curve closely matches the trajectory obtained via direct propagation in Cartesian coordinates. \cref{fig:r41error} illustrates the error in position and velocity from the two methods of propagation and confirms the accuracy of the high-fidelity M-GEqOE propagator. Finally, the Cartesian solution is transformed into the generalized coordinates, as demonstrated in \cref{fig:r41mgeqoe}. The consistency between the values of each orbital element over time is noted. Additionally, the periodic behavior of the M-GEqOEs aligns with the four orbital revolutions measured in the Earth-centered inertial frame, which together correspond to a single full revolution associated with the underlying sidereal resonance. A spacecraft in this orbit encounters the Earth four times with a perigee altitude of $92,415 \ km$. Each perigee pass is associated with an increased level of nonlinearity, which may introduce inaccuracies in the estimation of the spacecraft's states when restrictive uncertainty representations are employed.

\begin{figure}[!htbp]
    \centering
    \subfloat[HZ statistic]
        {\includegraphics[width=0.49\textwidth]{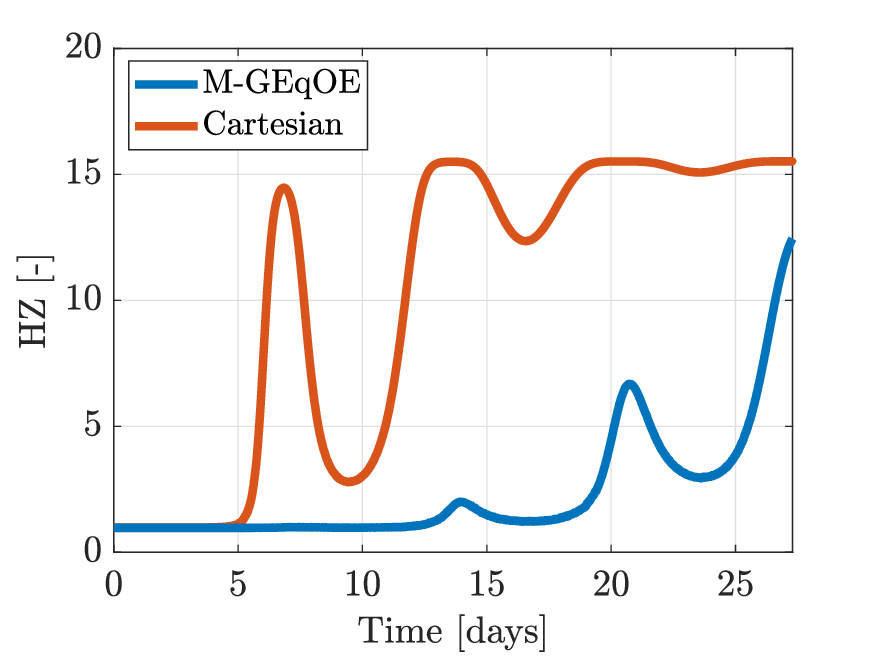}
        \label{fig:r41hz}}
    \subfloat[p-value]
        {\includegraphics[width=0.49\textwidth]{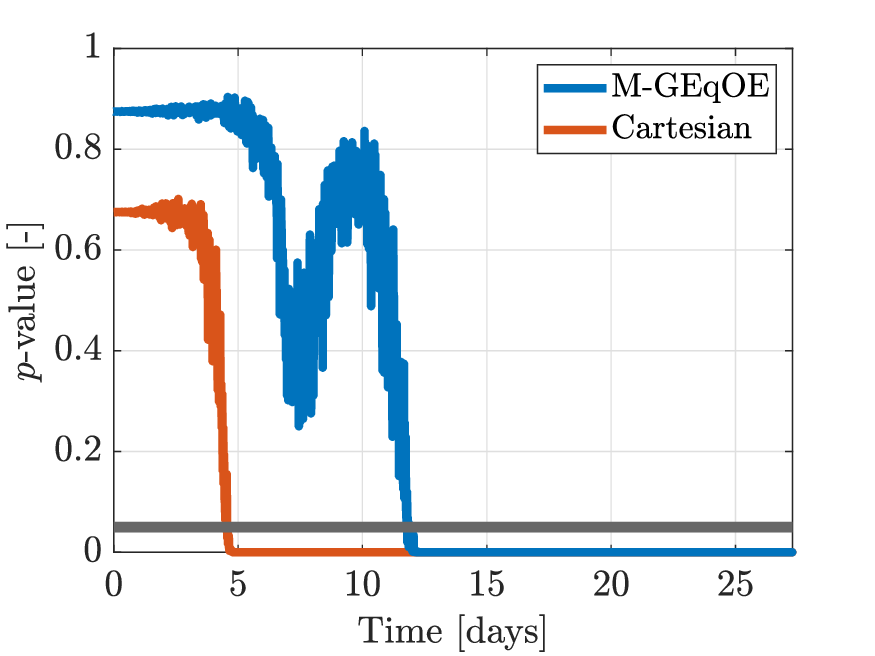}
        \label{fig:r41pval}}  
    \caption{Henze--Zirkler test applied to uncertainty propagated along the 4:1 resonant orbit.}
    \label{fig:r41unc}
\end{figure}

Uncertainty along the \rt{4}{1} resonant orbit is propagated in both generalized and Cartesian coordinates for one period, assuming initial $1\sigma$ values of uncertainty of $1 \ km$ and $1 \ cm/s$ along the position and velocity channels, respectively. The results for the HZ test appear in \cref{fig:r41unc}. For this particular orbit, propagation is initiated at perigee. In the Cartesian case, the departure from Gaussian behavior is evident at the next perigee crossing, which occurs approximately $6.8 \ days$ later. Further downstream, a drift in the value of the Cartesian HZ statistic is observed, and no locations along the orbit align with Gaussian behavior. For the M-GEqOE approach, Gaussianity is preserved up till the $14$-day mark, which corresponds to the third perigee crossing along the trajectory. Even in the M-GEqOE case, departures from purely Gaussian behavior appear at other locations along the orbit as well, though the extent of these deviations is significantly smaller than in the Cartesian approach. 

\begin{figure}[!htbp]
    \centering
    {
    \begin{tikzpicture}
        \node (popfig) at ([xshift=0cm,yshift=0cm] current page.center) {\includegraphics[width=1.05\textwidth]{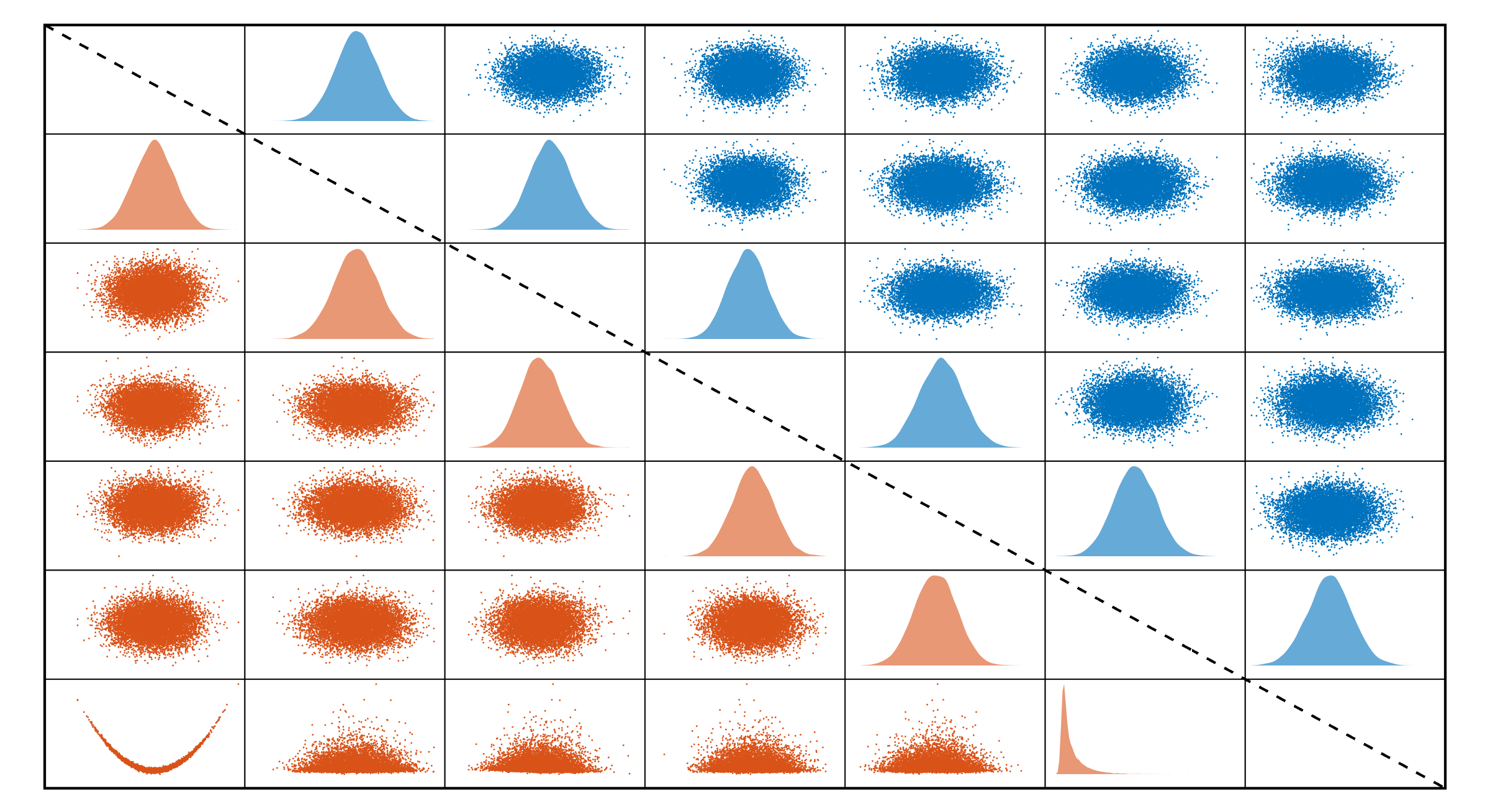}};
        \node at ([xshift=-6.75cm,yshift=1.0cm] popfig.center) {\large {$\lambda_{y}$}};
        \node at ([xshift=-6.75cm,yshift=0cm] popfig.center) {\large {$\lambda_{z}$}};
        \node at ([xshift=-6.75cm,yshift=-1.0cm] popfig.center) {\large {$\lambda_{\dot{x}}$}};
        \node at ([xshift=-6.75cm,yshift=-2.0cm] popfig.center) {\large {$\lambda_{\dot{y}}$}};
        \node at ([xshift=-6.75cm,yshift=-3.0cm] popfig.center) {\large {$\lambda_{\dot{z}}$}};
        \node at ([xshift=-5.4cm,yshift=-3.9cm] popfig.center) {\large {$\lambda_{x}$}};
        \node at ([xshift=-3.6cm,yshift=-3.9cm] popfig.center) {\large {$\lambda_{y}$}};
        \node at ([xshift=-1.8cm,yshift=-3.9cm] popfig.center) {\large {$\lambda_{z}$}};
        \node at ([xshift=0cm,yshift=-3.9cm] popfig.center) {\large {$\lambda_{\dot{x}}$}};
        \node at ([xshift=1.8cm,yshift=-3.9cm] popfig.center) {\large {$\lambda_{\dot{y}}$}};
        \node at ([xshift=3.6cm,yshift=-3.9cm] popfig.center) {\large {$\lambda_{\dot{z}}$}};
        \node at ([xshift=-3.6cm,yshift=3.9cm] popfig.center) {\large {$\lambda_{\tilde{p}}$}};
        \node at ([xshift=-1.8cm,yshift=3.9cm] popfig.center) {\large {$\lambda_{p_1}$}};
        \node at ([xshift=0cm,yshift=3.9cm] popfig.center) {\large {$\lambda_{p_2}$}};
        \node at ([xshift=1.8cm,yshift=3.9cm] popfig.center) {\large {$\lambda_{q_1}$}};
        \node at ([xshift=3.6cm,yshift=3.9cm] popfig.center) {\large {$\lambda_{q_2}$}};
        \node at ([xshift=5.4cm,yshift=3.9cm] popfig.center) {\large {$\lambda_{L}$}};
        \node at ([xshift=6.8cm,yshift=3cm] popfig.center) {\large {$\lambda_{\tilde{p}}$}};
        \node at ([xshift=6.8cm,yshift=2cm] popfig.center) {\large {$\lambda_{p_1}$}};
        \node at ([xshift=6.8cm,yshift=1cm] popfig.center) {\large {$\lambda_{p_2}$}};
        \node at ([xshift=6.8cm,yshift=0cm] popfig.center) {\large {$\lambda_{q_1}$}};
        \node at ([xshift=6.8cm,yshift=-1cm] popfig.center) {\large {$\lambda_{q_2}$}};
    \end{tikzpicture}}
    \caption{Eigenspace pairs plot for the 4:1 resonant orbit at ${t = 6.8 \ days}$ showing projections in Cartesian (lower triangular) and M-GEqOE (upper triangular) coordinates.}
    \label{fig:r41eigpop}
\end{figure}

\cref{fig:r41eigpop} presents the uncertainty clouds, transformed into the covariance eigenspace, for the M-GEqOE and Cartesian formulations at $t = 6.8 \ days$ into the propagation, corresponding to the second perigee pass. At this stage, the Cartesian uncertainty distribution shows the initial onset of non-Gaussian behavior, most visibly along the $\dot{z}$ component, as shown in orange in \cref{fig:r41eigpop}. While the underlying trajectory is nearly planar, a visible departure from Gaussian behavior is apparent along this direction under high-fidelity dynamics. On the other hand, the M-GEqOE ensemble plotted in blue preserves Gaussianity, consistent with the results of the HZ test.

\subsection{Elliptical Lunar Frozen Orbit}
\begin{figure}[!htbp]
    \centering
    \subfloat[Earth-Moon rotating frame]
        {\includegraphics[width=0.48\textwidth]{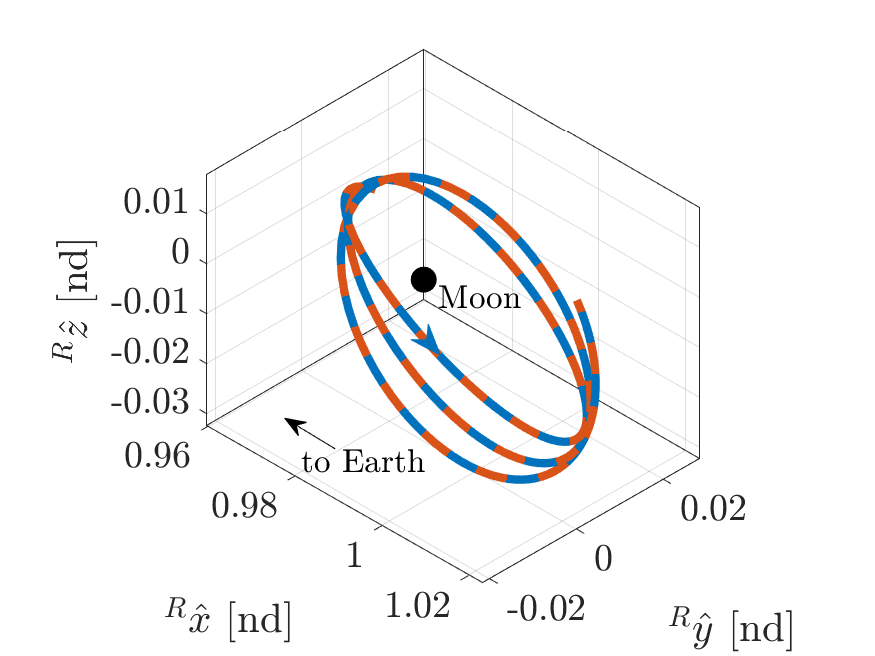}
        \label{fig:elforot}}
    \subfloat[Moon-centered inertial frame]
        {\includegraphics[width=0.48\textwidth]{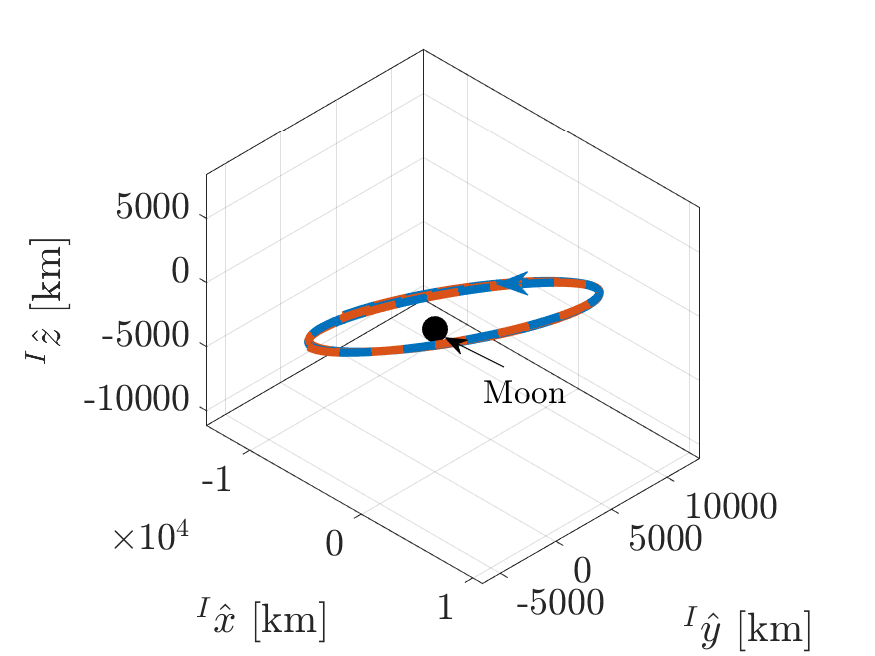}
        \label{fig:elfomci}}   \\
    \subfloat[Position and velocity errors]
        {\includegraphics[width=0.6\textwidth]{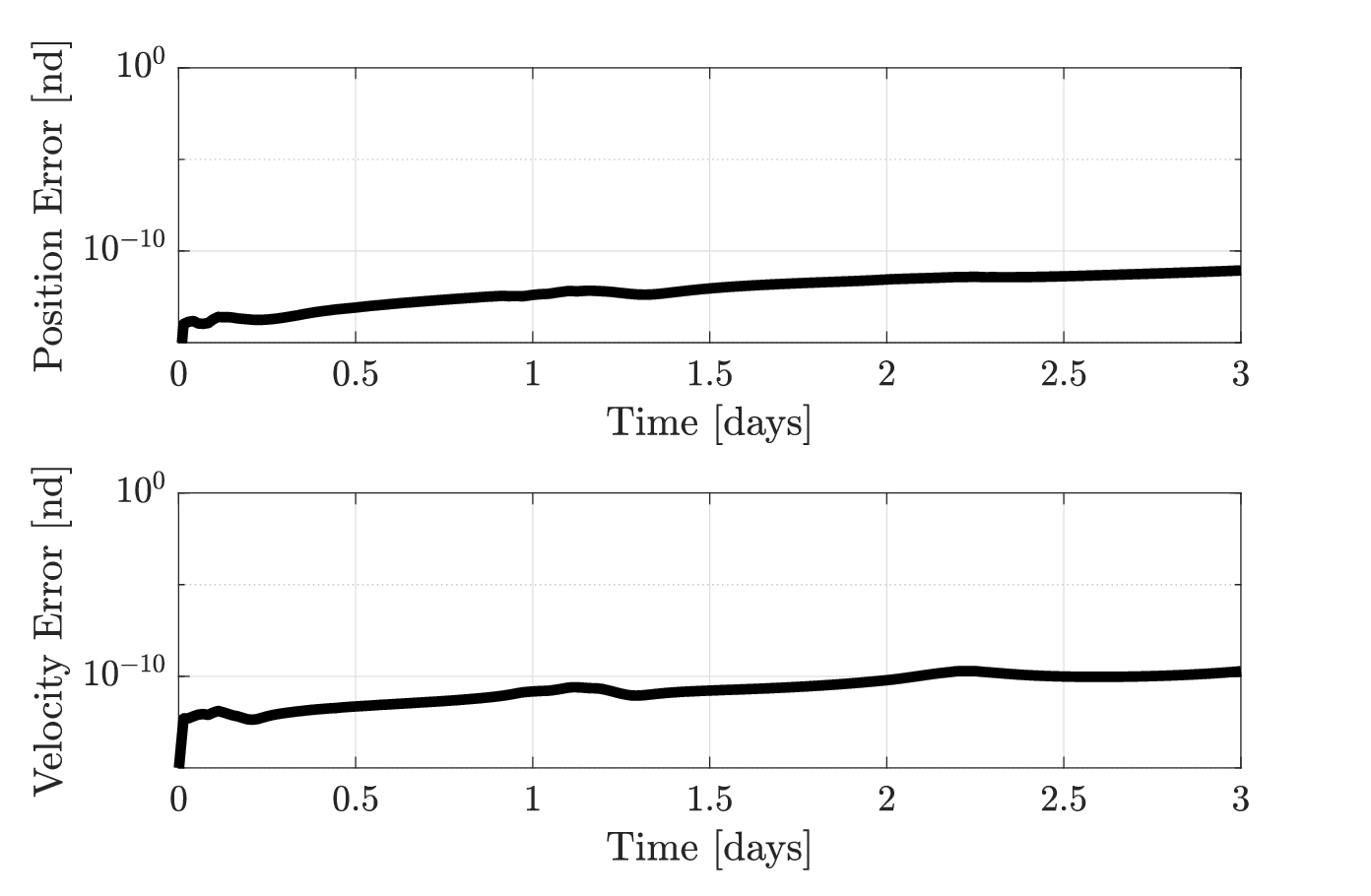}
        \label{fig:elfoerror}}   \\    
    \subfloat[Evolution of M-GEqOEs]
    {\includegraphics[width=0.9\textwidth]{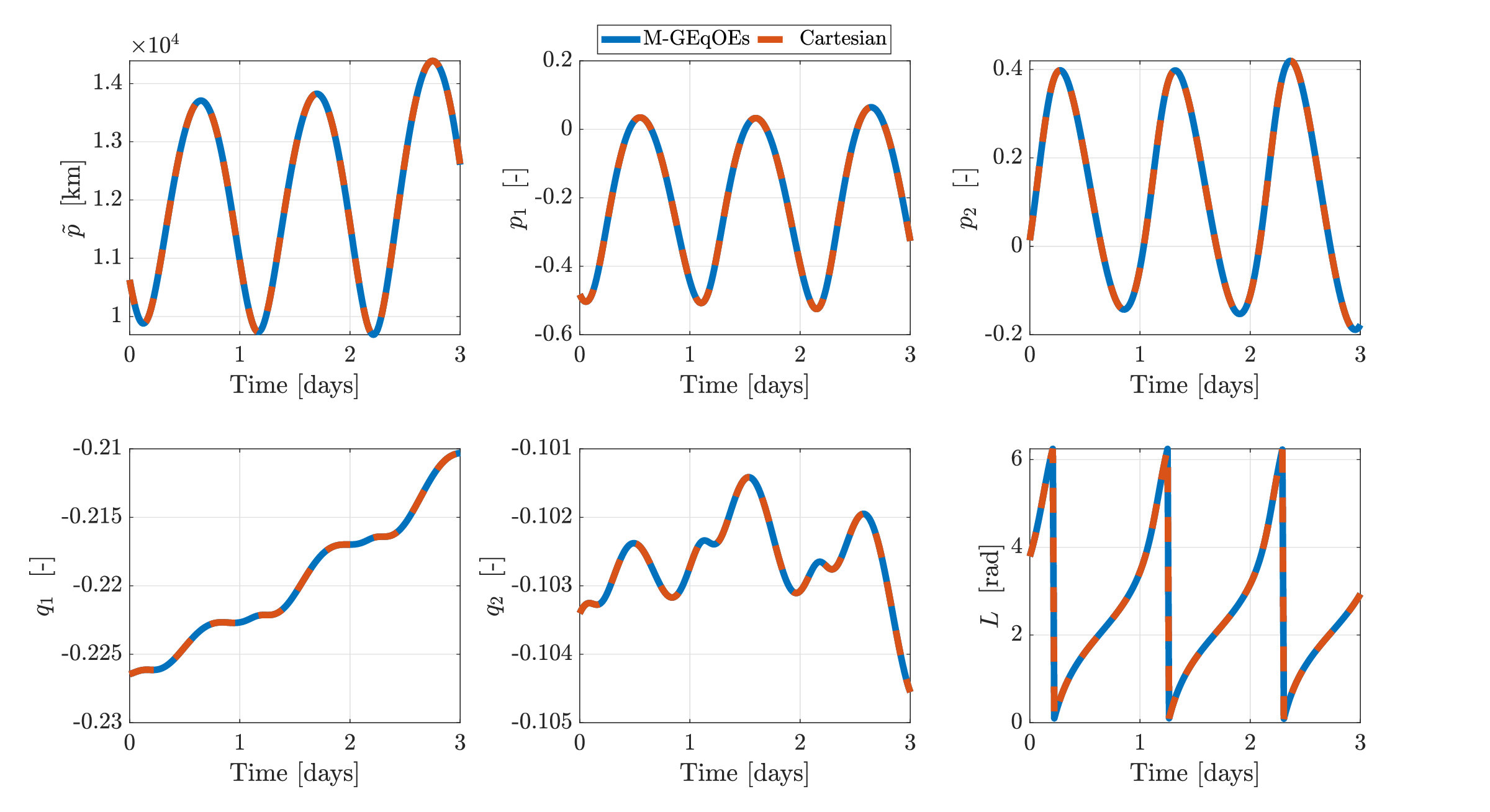}\label{fig:elfomgeqoe}}
    \caption{Cartesian and M-GEqOE representations of the ELFO over $3 \ days$. Blue and orange curves represent the M-GEqOE and Cartesian solutions, respectively.}
    \label{fig:elfostate}
\end{figure}

As a third example for high-fidelity state and uncertainty propagation in the cislunar domain, consider an Elliptic Lunar Frozen Orbit (ELFO) for operations in the lunar vicinity, including potentially hosting data relay constellations \citep{Ely2006}. These orbits are characterized by high values of inclination, eccentricity, and altitude, allowing long-term coverage of the lunar south pole. More recently, extensive research has been done to refine these orbits for rapid constellation design in a higher-fidelity dynamical environment \citep{Singh2020,Park2025}. In the current work, the initial conditions of a representative ELFO, as provided by \cite{Park2025}, are used for analysis. The initial state is transformed into the generalized coordinates for high-fidelity propagation over $30 \ days$ using the M-GEqOE equations. The resulting solution, along with the trajectory obtained via direct Cartesian propagation, appears in \cref{fig:elforot,fig:elfomci} as viewed in the Earth-Moon rotating and Moon-centered inertial frames. The accuracy of the M-GEqOE solution is confirmed via \cref{fig:elfoerror}, which illustrates the error in position and velocity over $3 \ days$. \cref{fig:elfomgeqoe} illustrates the M-GEqOEs that characterize the ELFO, where the orange curves represent the Cartesian-propagated states transformed into the generalized coordinates. 

\begin{figure}[!htbp]
    \centering
    \subfloat[HZ statistic]
        {\includegraphics[width=0.49\textwidth]{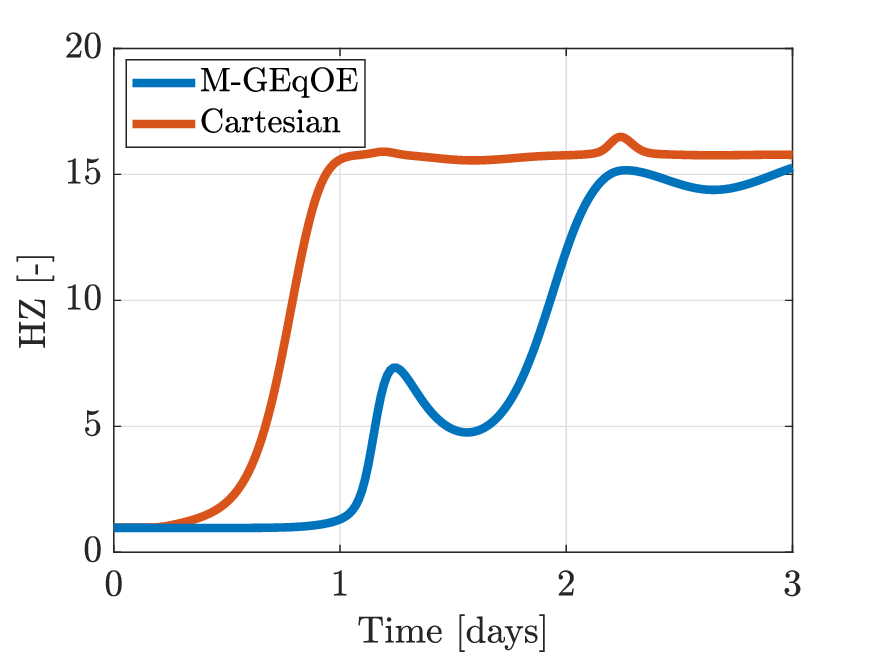}
        \label{fig:elfohz}}
    \subfloat[p-value]
        {\includegraphics[width=0.49\textwidth]{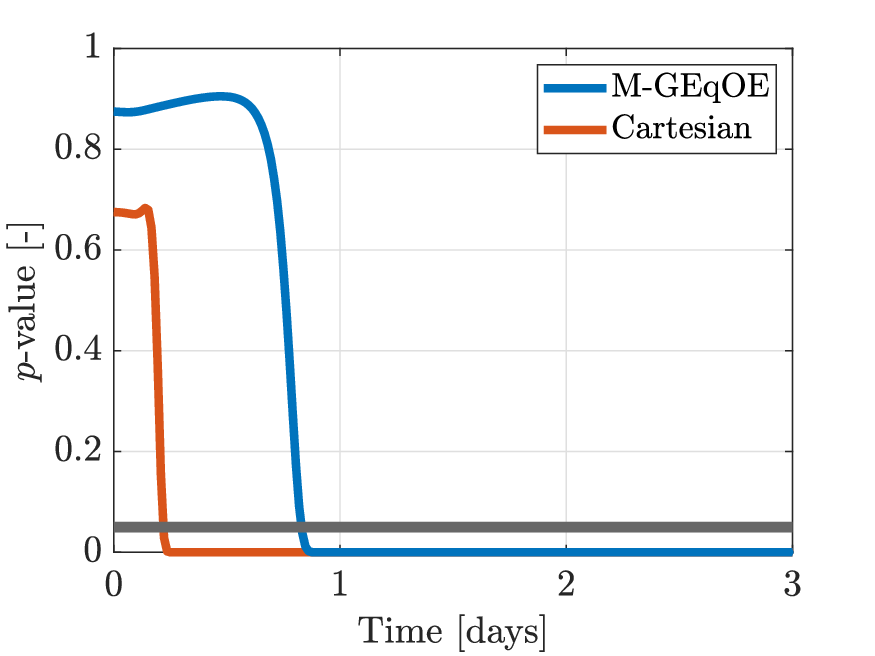}
        \label{fig:elfopval}}  
    \caption{Henze--Zirkler test applied to uncertainty propagated along the ELFO.}
    \label{fig:elfounc}
\end{figure}

To assess the evolution of uncertainty along the selected ELFO, initial $1\sigma$ uncertainty values of $1 \ km$ and $1 \ cm/s$ are assumed along the position and velocity components. The HZ statistic and the associated $p$-value over $3 \ days$ of propagation appear in \cref{fig:elfounc}. The M-GEqOE methodology is able to maintain Gaussianity for approximately $1 \ day$, whereas the Cartesian representation departs from Gaussian behavior much earlier. In the Cartesian case, the HZ statistic increases sharply and the associated $p$-value rapidly collapses to zero, remaining there for the remainder of the propagation. Beyond this point, non-Gaussianity is apparent in both coordinate sets, though the deviations are less severe in the generalized coordinates compared to the Cartesian case for approximately $2 \ days$ into the propagation.

\begin{figure}[!htbp]
    \centering
    {
    \begin{tikzpicture}
        \node (popfig) at ([xshift=0cm,yshift=0cm] current page.center) {\includegraphics[width=1.05\textwidth]{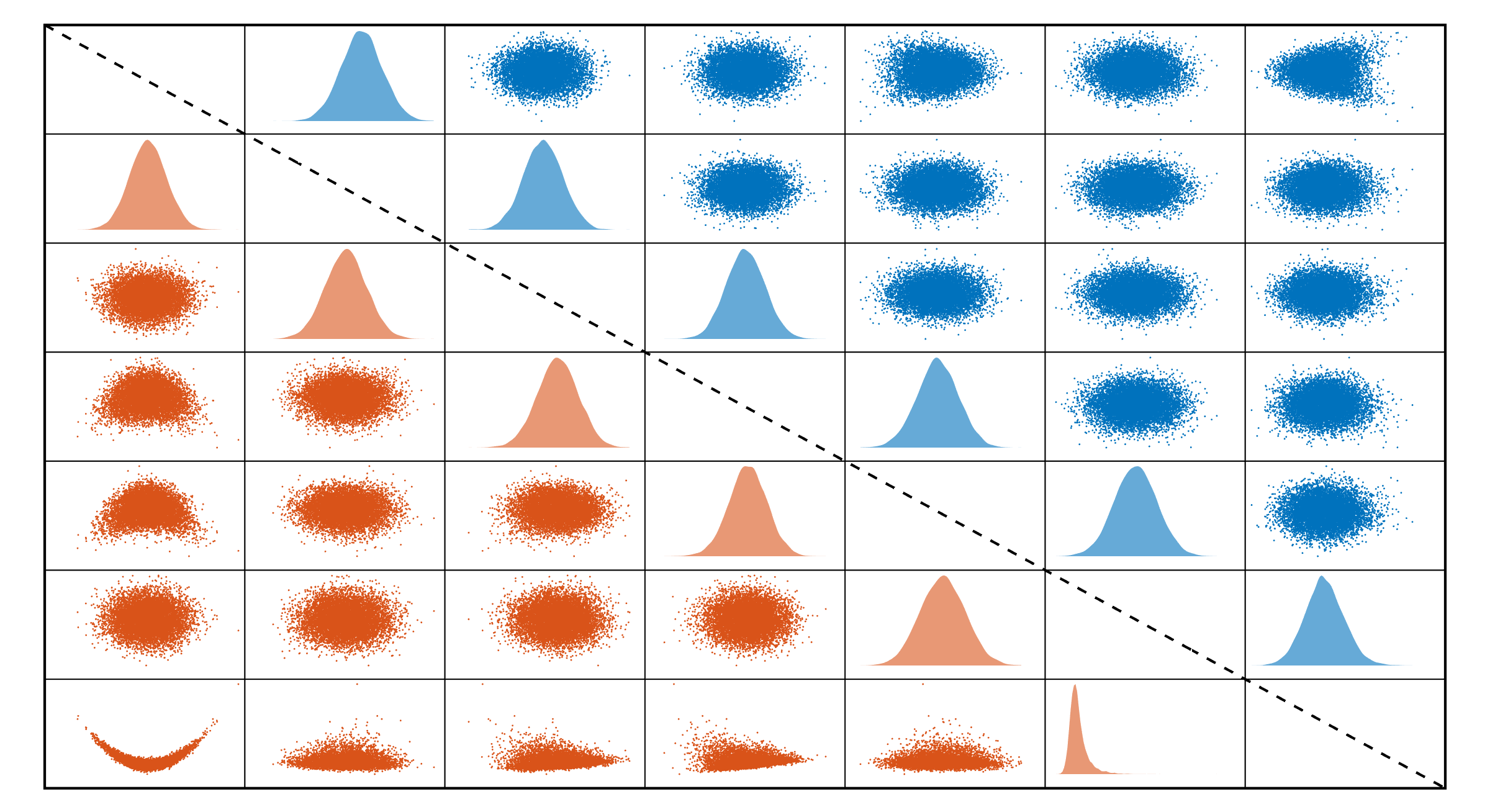}};
        \node at ([xshift=-6.75cm,yshift=1.0cm] popfig.center) {\large {$\lambda_{y}$}};
        \node at ([xshift=-6.75cm,yshift=0cm] popfig.center) {\large {$\lambda_{z}$}};
        \node at ([xshift=-6.75cm,yshift=-1.0cm] popfig.center) {\large {$\lambda_{\dot{x}}$}};
        \node at ([xshift=-6.75cm,yshift=-2.0cm] popfig.center) {\large {$\lambda_{\dot{y}}$}};
        \node at ([xshift=-6.75cm,yshift=-3.0cm] popfig.center) {\large {$\lambda_{\dot{z}}$}};
        \node at ([xshift=-5.4cm,yshift=-3.9cm] popfig.center) {\large {$\lambda_{x}$}};
        \node at ([xshift=-3.6cm,yshift=-3.9cm] popfig.center) {\large {$\lambda_{y}$}};
        \node at ([xshift=-1.8cm,yshift=-3.9cm] popfig.center) {\large {$\lambda_{z}$}};
        \node at ([xshift=0cm,yshift=-3.9cm] popfig.center) {\large {$\lambda_{\dot{x}}$}};
        \node at ([xshift=1.8cm,yshift=-3.9cm] popfig.center) {\large {$\lambda_{\dot{y}}$}};
        \node at ([xshift=3.6cm,yshift=-3.9cm] popfig.center) {\large {$\lambda_{\dot{z}}$}};
        \node at ([xshift=-3.6cm,yshift=3.9cm] popfig.center) {\large {$\lambda_{\tilde{p}}$}};
        \node at ([xshift=-1.8cm,yshift=3.9cm] popfig.center) {\large {$\lambda_{p_1}$}};
        \node at ([xshift=0cm,yshift=3.9cm] popfig.center) {\large {$\lambda_{p_2}$}};
        \node at ([xshift=1.8cm,yshift=3.9cm] popfig.center) {\large {$\lambda_{q_1}$}};
        \node at ([xshift=3.6cm,yshift=3.9cm] popfig.center) {\large {$\lambda_{q_2}$}};
        \node at ([xshift=5.4cm,yshift=3.9cm] popfig.center) {\large {$\lambda_{L}$}};
        \node at ([xshift=6.8cm,yshift=3cm] popfig.center) {\large {$\lambda_{\tilde{p}}$}};
        \node at ([xshift=6.8cm,yshift=2cm] popfig.center) {\large {$\lambda_{p_1}$}};
        \node at ([xshift=6.8cm,yshift=1cm] popfig.center) {\large {$\lambda_{p_2}$}};
        \node at ([xshift=6.8cm,yshift=0cm] popfig.center) {\large {$\lambda_{q_1}$}};
        \node at ([xshift=6.8cm,yshift=-1cm] popfig.center) {\large {$\lambda_{q_2}$}};
    \end{tikzpicture}}
    \caption{Eigenspace pairs plot for the ELFO at ${t = 1 \ day}$ showing projections in Cartesian (lower triangular) and M-GEqOE (upper triangular) coordinates.}
    \label{fig:elfoeigpop}
\end{figure}
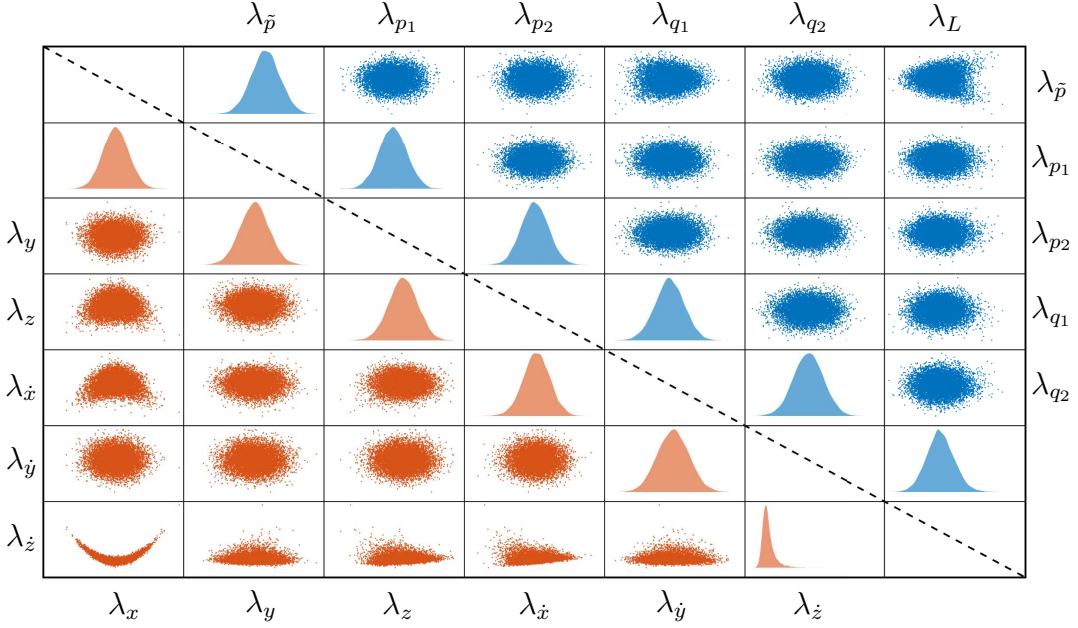

\cref{fig:elfoeigpop} illustrates the eigenspace uncertainty clouds for uncertainty propagated in the M-GEqOE (blue) and Cartesian (orange) coordinates at $1 \ day$ downstream. In the M-GEqOE coordinates, the particle cloud exhibits Gaussian structure across most projections, with some curvature appearing in a few components as the nonlinear effects accumulate. On the other hand, the Cartesian coordinates exhibit elongation and curvature across multiple components, consistent with stronger departure from Gaussianity at this stage of the propagation. 

\subsection{3:1 Spatial Sidereal Resonant Orbit}
\begin{figure}[!htbp]
    \centering
    \subfloat[Earth-Moon rotating frame]
        {\includegraphics[width=0.48\textwidth]{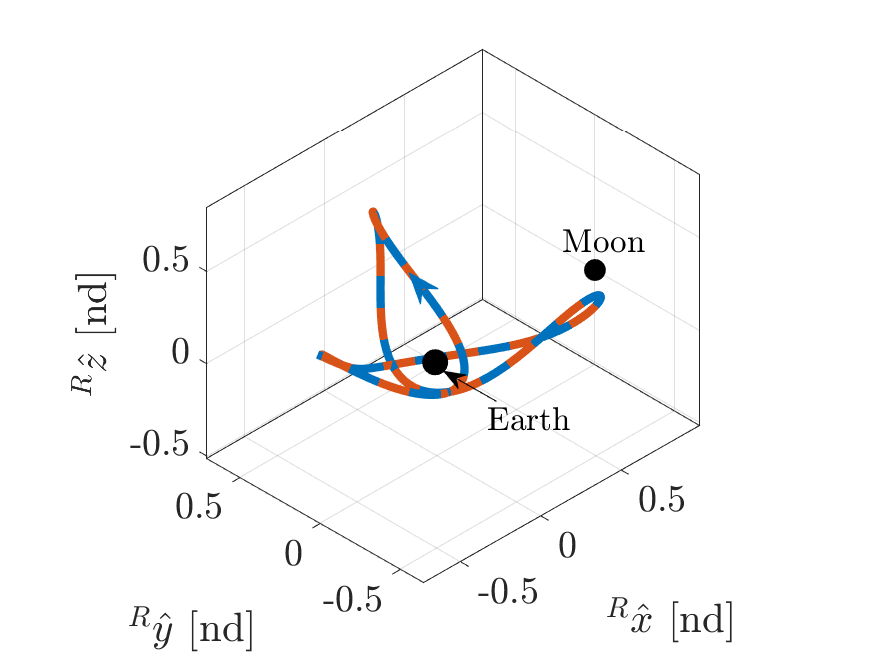}
        \label{fig:r31rot}}
    \subfloat[Earth-centered inertial frame]
        {\includegraphics[width=0.48\textwidth]{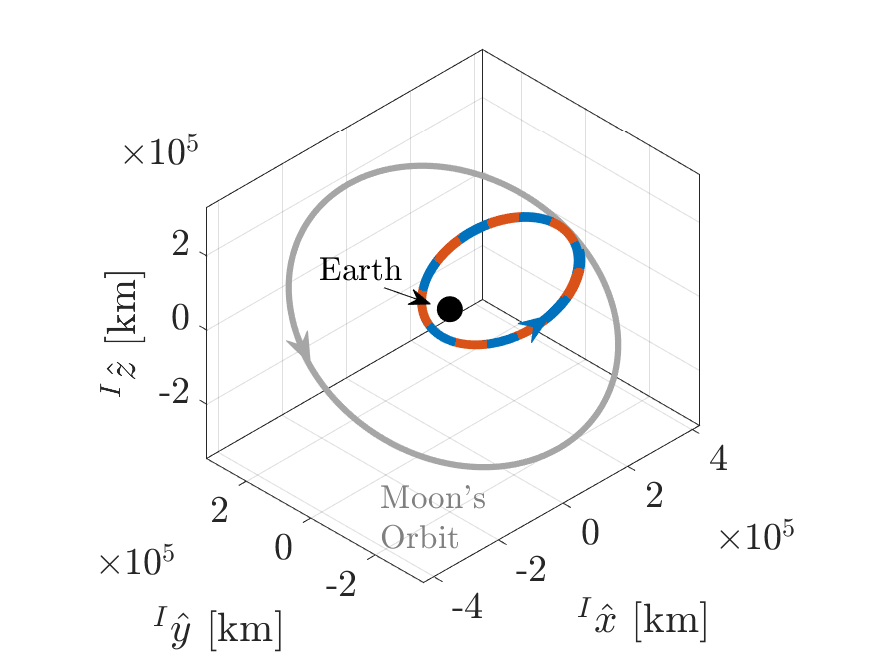}
        \label{fig:r31eci}}   \\
    \subfloat[Position and velocity errors]
        {\includegraphics[width=0.6\textwidth]{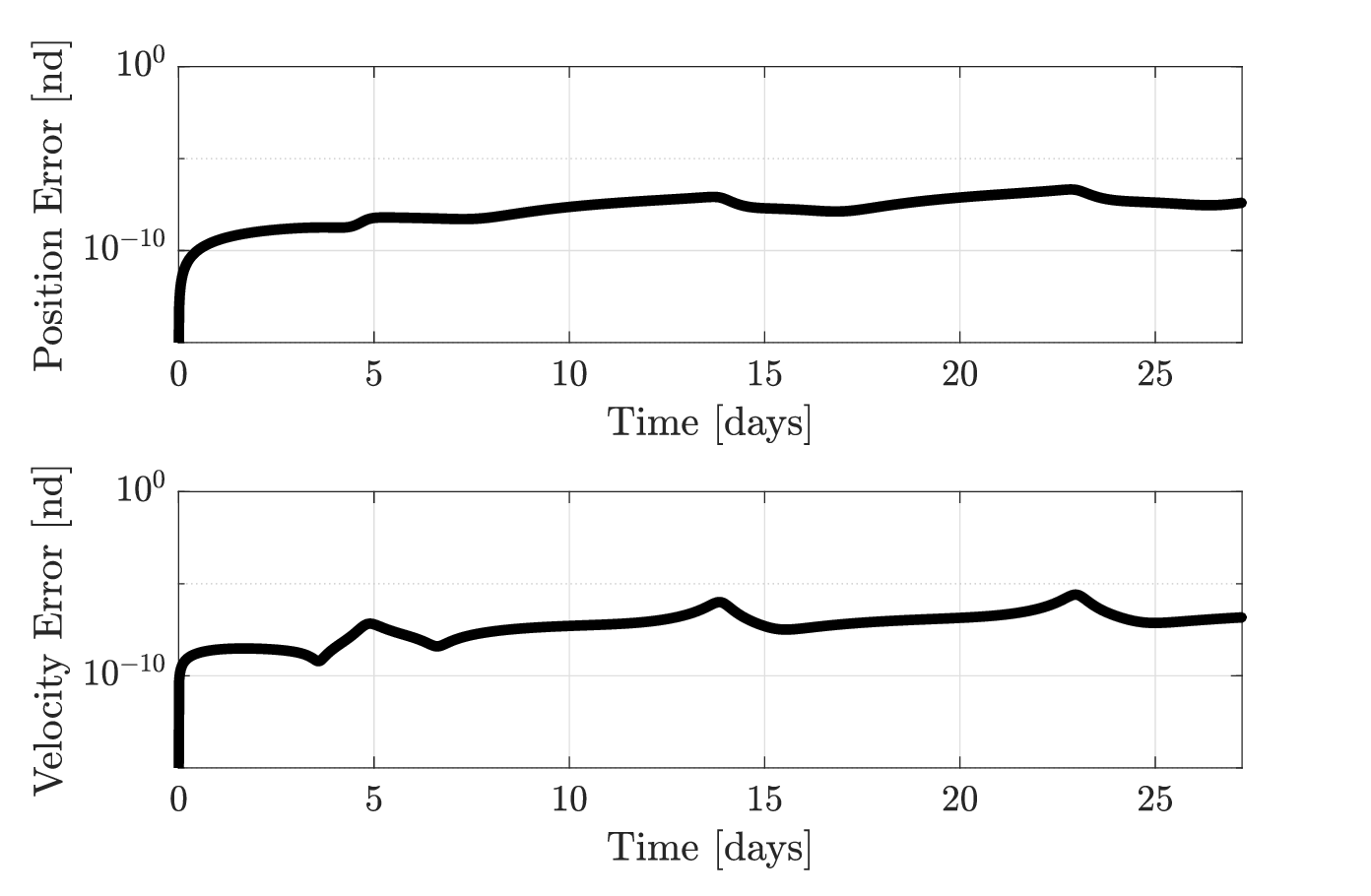}
        \label{fig:r31error}}   \\    
    \subfloat[Evolution of M-GEqOEs]
    {\includegraphics[width=0.9\textwidth]{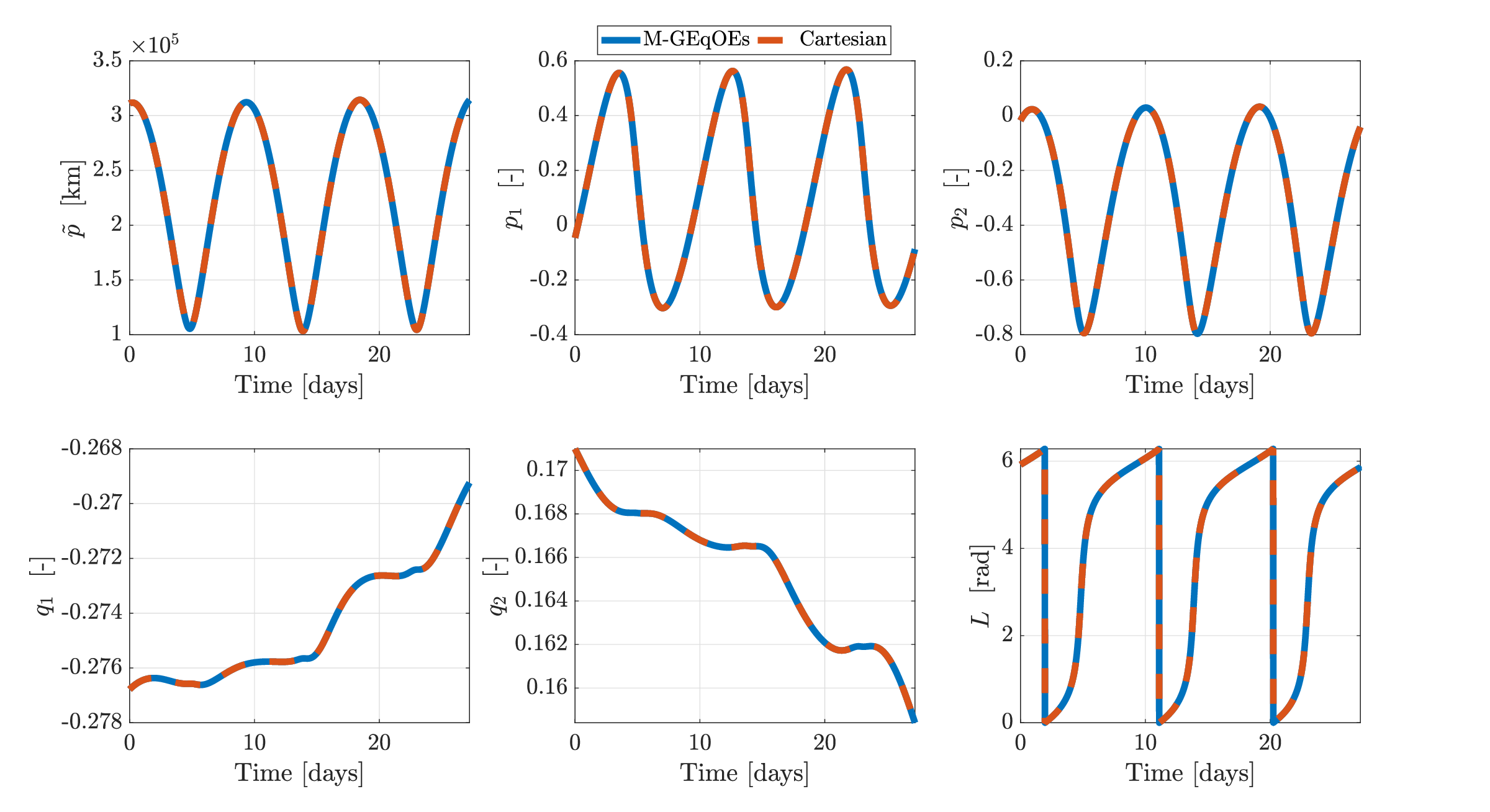}\label{fig:r31mgeqoe}}
    \caption{Cartesian and M-GEqOE representations of the 3:1 sidereal resonant orbit over $27.22\ days$. Blue and orange curves represent the M-GEqOE and Cartesian solutions, respectively.}
    \label{fig:r31state}
\end{figure}

As a final representative example, another Earth-centered sidereal resonant orbit from the family of spatial \rt{3}{1} resonant orbits is selected. This particular orbit possesses a geometry similar to that of the Interstellar Boundary Explorer (IBEX) extended mission \cite{Dichmann2013}. The period of the selected orbit is approximately $27.22 \ days$. \cref{fig:r31rot,fig:r31eci} illustrate the trajectories obtained via the M-GEqOE propagation (blue) and the Cartesian propagation (orange) in the Earth-Moon rotating and Earth-centered inertial frames. The position and velocity errors between the two propagation methodologies, shown in \cref{fig:r31error}, initially agree with the levels observed in the previous examples, followed by an order of magnitude increase associated with passing through perigee. \cref{fig:r31mgeqoe} illustrates the time history of the M-GEqOEs for this orbit, where the periodicity of the elements is consistent with the sidereal resonance of this orbit. 

\begin{figure}[!htbp]
    \centering
    \subfloat[HZ statistic]
        {\includegraphics[width=0.4925\textwidth]{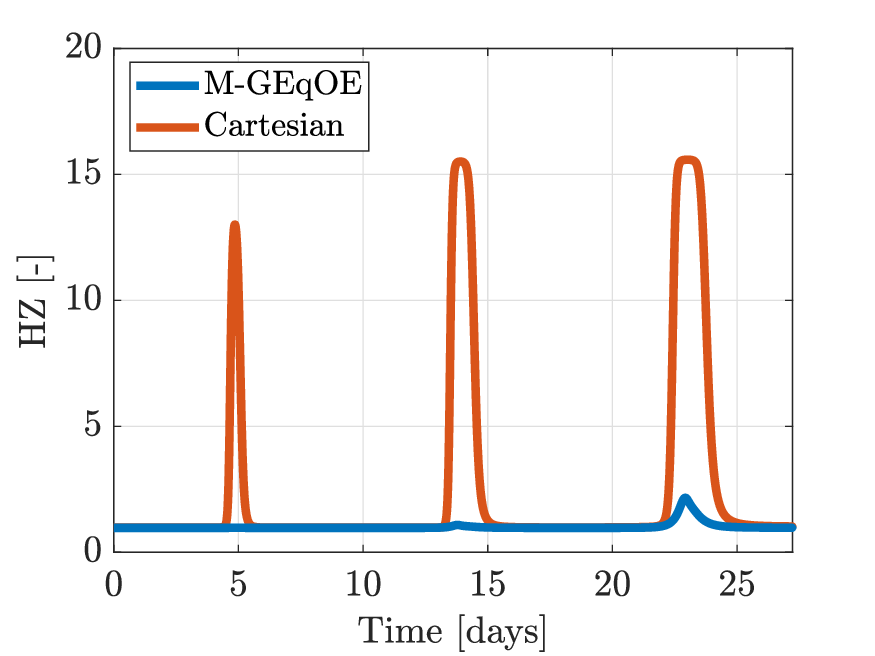}
        \label{fig:r31hz}}
    \subfloat[p-value]
        {\includegraphics[width=0.4925\textwidth]{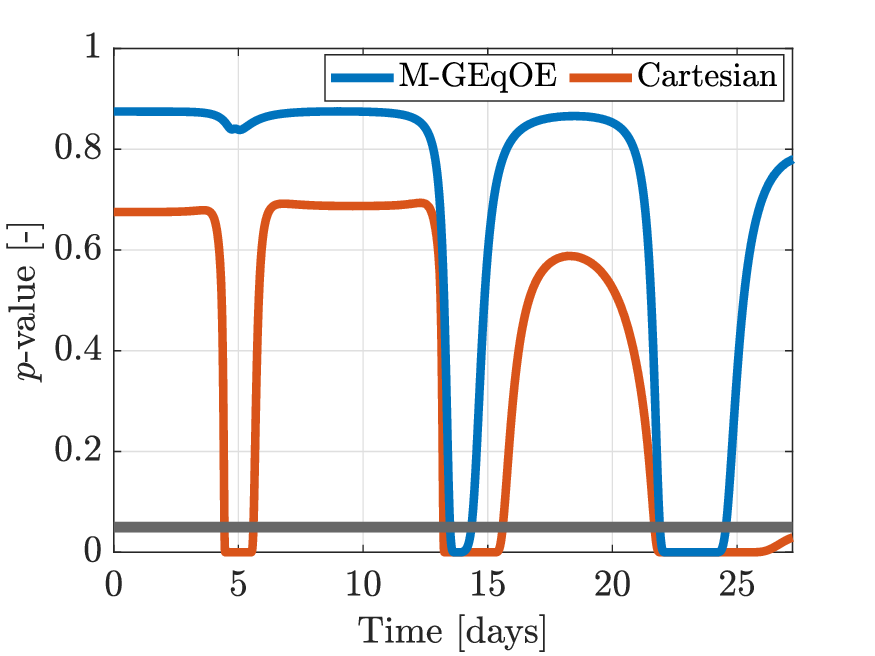}
        \label{fig:r31pval}}  
    \caption{Henze--Zirkler test applied to uncertainty propagated along the 3:1 resonant orbit.}
    \label{fig:r31unc}
\end{figure}

Uncertainty is propagated for the \rt{3}{1} resonant orbit, assuming $1\sigma$ standard deviations of $10 \ km$ and $10 \ cm/s$ along the position and velocity channels, respectively. \cref{fig:r31unc} illustrates the HZ statistic and the associated $p$-value computed along the 3:1 resonant orbit for both M-GEqOE (blue) and Cartesian (orange) representations. Uncertainty propagated in Cartesian coordinates exhibits pronounced non-Gaussianity at each perigee pass, as indicated by sharp increases in the HZ statistic and corresponding drops in the $p$-value below the selected significance level. On the other hand, the M-GEqOE representation is able to maintain the nominal Gaussian behavior, with $p$-values consistently above the Gaussian rejection threshold over most of the trajectory. Although non-Gaussianity is observed at the second and third perigee passes in the M-GEqOE coordinates, the HZ statistic exhibits a smaller increase and recovers more rapidly compared to the Cartesian case, demonstrating that the M-GEqOE methodology better preserves Gaussianity over time. 

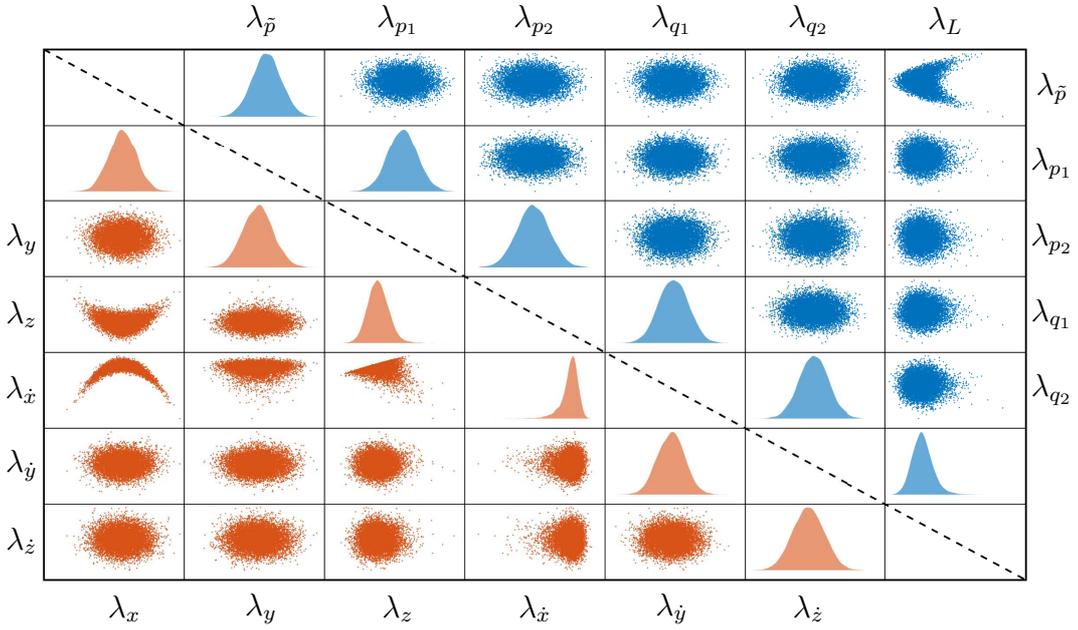
\begin{figure}[!htbp]
    \centering
    {
    \begin{tikzpicture}
        \node (popfig) at ([xshift=0cm,yshift=0cm] current page.center) {\includegraphics[width=1.05\textwidth]{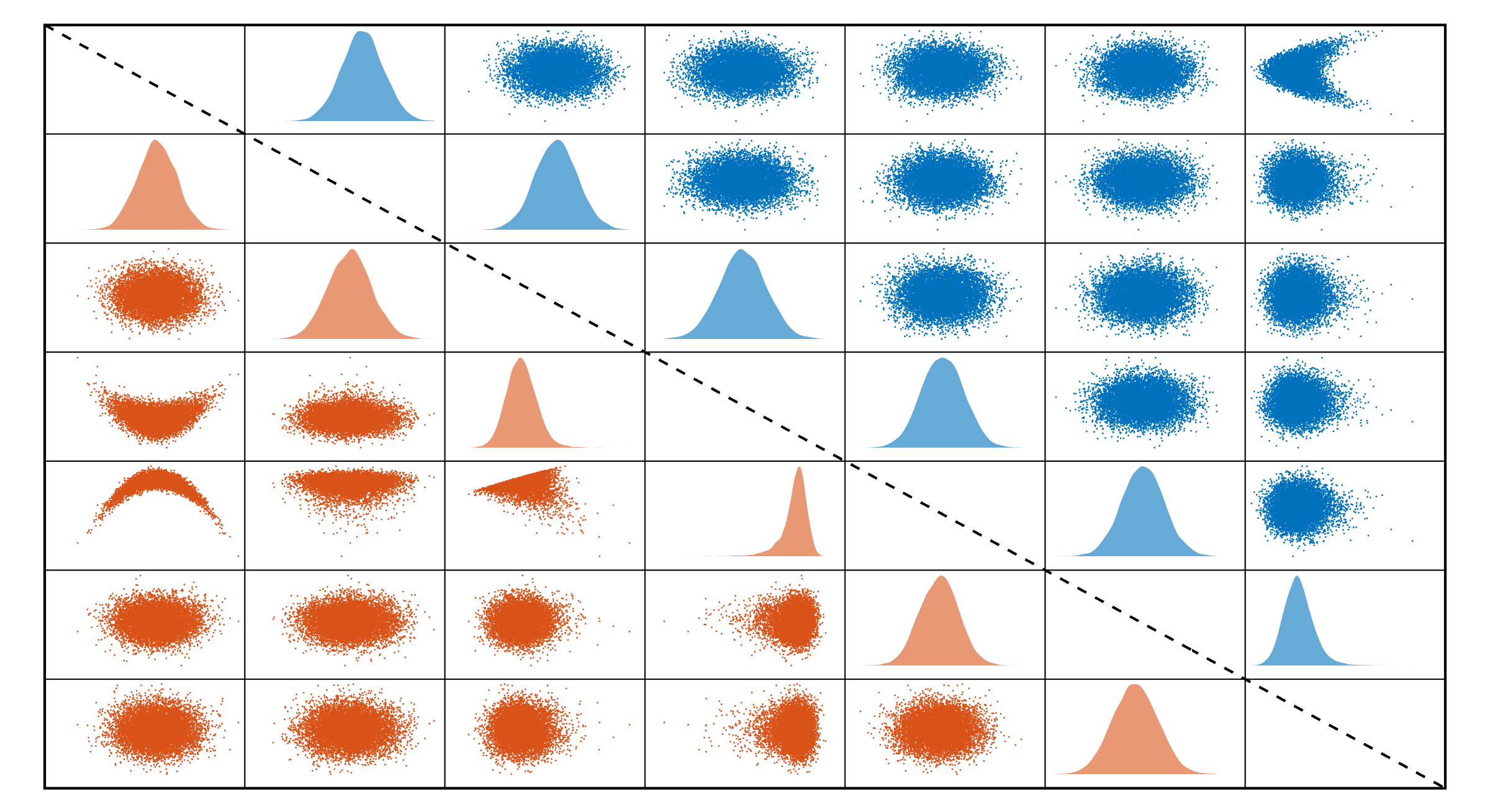}};
        \node at ([xshift=-6.75cm,yshift=1.0cm] popfig.center) {\large {$\lambda_{y}$}};
        \node at ([xshift=-6.75cm,yshift=0cm] popfig.center) {\large {$\lambda_{z}$}};
        \node at ([xshift=-6.75cm,yshift=-1.0cm] popfig.center) {\large {$\lambda_{\dot{x}}$}};
        \node at ([xshift=-6.75cm,yshift=-2.0cm] popfig.center) {\large {$\lambda_{\dot{y}}$}};
        \node at ([xshift=-6.75cm,yshift=-3.0cm] popfig.center) {\large {$\lambda_{\dot{z}}$}};
        \node at ([xshift=-5.4cm,yshift=-3.9cm] popfig.center) {\large {$\lambda_{x}$}};
        \node at ([xshift=-3.6cm,yshift=-3.9cm] popfig.center) {\large {$\lambda_{y}$}};
        \node at ([xshift=-1.8cm,yshift=-3.9cm] popfig.center) {\large {$\lambda_{z}$}};
        \node at ([xshift=0cm,yshift=-3.9cm] popfig.center) {\large {$\lambda_{\dot{x}}$}};
        \node at ([xshift=1.8cm,yshift=-3.9cm] popfig.center) {\large {$\lambda_{\dot{y}}$}};
        \node at ([xshift=3.6cm,yshift=-3.9cm] popfig.center) {\large {$\lambda_{\dot{z}}$}};
        \node at ([xshift=-3.6cm,yshift=3.9cm] popfig.center) {\large {$\lambda_{\tilde{p}}$}};
        \node at ([xshift=-1.8cm,yshift=3.9cm] popfig.center) {\large {$\lambda_{p_1}$}};
        \node at ([xshift=0cm,yshift=3.9cm] popfig.center) {\large {$\lambda_{p_2}$}};
        \node at ([xshift=1.8cm,yshift=3.9cm] popfig.center) {\large {$\lambda_{q_1}$}};
        \node at ([xshift=3.6cm,yshift=3.9cm] popfig.center) {\large {$\lambda_{q_2}$}};
        \node at ([xshift=5.4cm,yshift=3.9cm] popfig.center) {\large {$\lambda_{L}$}};
        \node at ([xshift=6.8cm,yshift=3cm] popfig.center) {\large {$\lambda_{\tilde{p}}$}};
        \node at ([xshift=6.8cm,yshift=2cm] popfig.center) {\large {$\lambda_{p_1}$}};
        \node at ([xshift=6.8cm,yshift=1cm] popfig.center) {\large {$\lambda_{p_2}$}};
        \node at ([xshift=6.8cm,yshift=0cm] popfig.center) {\large {$\lambda_{q_1}$}};
        \node at ([xshift=6.8cm,yshift=-1cm] popfig.center) {\large {$\lambda_{q_2}$}};
    \end{tikzpicture}}
    \caption{Eigenspace pairs plot for the 3:1 resonant orbit at ${t = 23 \ days}$ showing projections in Cartesian (lower triangular) and M-GEqOE (upper triangular) coordinates.}
    \label{fig:r31eigpop}
\end{figure}

The eigenspace pairs plot in \cref{fig:r31eigpop}, constructed at $t = 23$ days and corresponding to the third perigee pass, further confirms these trends. Both representations exhibit non-Gaussian behavior at this point, but the deviations are significantly lower in the M-GEqOE set than in Cartesian coordinates. Thus, the M-GEqOE set is demonstrated to preserve Gaussianity for longer durations than Cartesian propagation of uncertainty, mitigating the effects of nonlinearities at perigee and providing an improvement in uncertainty realism compared to the conventional Cartesian approach. 

\section{Conclusion}
Cislunar space poses a challenging dynamical environment that requires high-fidelity dynamical modeling to accurately predict spacecraft trajectories and associated uncertainties. Space Domain Awareness in this regime relies heavily on accurately modeling the uncertainty in a spacecraft’s state and its evolution over time. Traditional assumptions of Gaussianity for the initial uncertainty and its downstream propagation are often overly restrictive and may fail to capture the true behavior of the system. As such, alternative methodologies are required to retain Gaussian representations while ensuring a realistic representation of uncertainty. These techniques are particularly critical for establishing and maintaining custody of both known and unidentified objects in the cislunar environment.

In the current work, the Modified Generalized Equinoctial Orbital Elements (M-GEqOEs) are explored as an option for cislunar state and uncertainty propagation under the presence of the necessary gravitational perturbations. The M-GEqOE set allows the direct inclusion of conservative perturbations, offering a low-complexity methodology for modeling dynamics in this complex regime. Various cislunar orbits are constructed to demonstrate the applicability of this element set for high-fidelity propagation under the gravitational influence of the Earth, Moon, and Sun. 

In addition to state propagation, uncertainty characterization via the M-GEqOE set is explored for various cislunar orbits. The Henze--Zirkler test for multivariate normality is employed to characterize the uncertainty over time in both M-GEqOE and Cartesian coordinates. For the orbits considered in this work, improved preservation of Gaussianity is observed via the M-GEqOE methodology, compared to the Cartesian approach. In particular, the generalized coordinate representation maintains Gaussian behavior over longer propagation intervals and exhibits improved uncertainty characterization in sensitive and highly nonlinear regions of the trajectories.

\subsubsection*{Acknowledgements} This work was supported by the Air Force Office of Scientific Research (AFOSR) under agreement number FA9550-23-1-0646, Create the Future Independent Research Effort (CFIRE).


\bibliography{bibliography}

\end{document}